\documentclass[11pt]{siamltex}

\usepackage[english]{babel}
\usepackage[dvips]{graphicx}
\usepackage{amssymb}
\usepackage{bm}
\usepackage{amsmath}
\usepackage{xcolor}
\usepackage{url}
\usepackage{hyperref}

\newcommand{\RR}{\mathbb{R}}

\newcommand{\PP}{\mathbb{P}}
\def\C{\mathcal{C}}

\newcommand{\la}{\lambda}

\newcommand{\vphi}{\varphi}
\newcommand{\ds}{\displaystyle}
\def\C{\mathcal{C}}
\def\PP{\vbox {\hbox{I\hskip-2.1pt P\hfil}}}
\def\RR{\vbox {\hbox{I\hskip-2.1pt R\hfil}}}

\newtheorem{remark}{Remark}
\newtheorem{example}{Example}

\begin{document}

\title{A global approximation method for second-kind nonlinear integral equations
}


\author{Luisa Fermo         \and
        Anna Lucia Laguardia \and Concetta Laurita \and
        Maria Grazia Russo
}


\author{L. Fermo \thanks{Department of Mathematics and Computer Science, University of Cagliari,
Via Ospedale 72, 09124 Cagliari, Italy, email: fermo@unica.it. }      
           \and
           A. L. Laguardia \thanks{Department of Mathematics and Computer Science, University of Basilicata, Viale dell'Ateneo Lucano 10, 85100 Potenza, Italy, email: annalucia.laguardia@unibas.it, concetta.laurita@unibas.it, mariagrazia.russo@unibas.it}  
               \and
               C. Laurita  \footnotemark[2]          
               \and
               M.G. Russo \footnotemark[2]
}

\date{Received: date / Accepted: date}

\maketitle

\begin{abstract}
A global approximation method of Nystr\"om type is explored for the numerical solution of a class of nonlinear integral equations of the second kind.
The cases of smooth and weakly singular kernels are both considered. In the first occurrence, the method uses a Gauss-Legendre rule whereas in the second one resorts to a product rule based on Legendre nodes.  Stability and convergence are proved in functional spaces equipped with the uniform norm and several numerical tests are given to show the good performance of the proposed method. An application to the interior Neumann problem for the Laplace equation with nonlinear boundary conditions is also considered.

\begin{keywords}
Nonlinear integral equations, Hammerstein-type integral equations, \\Nystr\"om method, Boundary integral equations.
\end{keywords} 
\begin{AMS}
65R20, 45G10, 47H30
\end{AMS}
\end{abstract}

\section{Introduction}
The goal of this paper is to develop a numerical method for the following integral equation
\begin{equation}
\label{EqInFr}
f(y) -  \int_{-1}^{1} k_1(x,y)f(x) dx-  \int_{-1}^{1} k_2(x,y)h(x,f(x)) dx= g(y), \quad y \in [-1,1],
\end{equation}
where $f$ is to be determined, $k_1$, $k_2$ and $g$ are  given functions, and $h(x,v)$ is  a known function which is assumed to be nonlinear in $v$.

Integral equations of type \eqref{EqInFr} have wide applications in models involving nonlinearities such as heat radiation, heat transfer, acoustic, elasticity, and electromagnetic problems; see \cite{bialecki1981,kelmanson1984}.
Some of these models are mathematically represented in terms of boundary value problems having nonlinear boundary conditions which can be reformulated in terms of \eqref{EqInFr}; see \cite{ruotsalainen1988} and Section \ref{sec:BIE}. In applicative contexts, the kernels $k_1$ and $k_2$ of \eqref{EqInFr} are smooth and/or weakly singular. For instance, in the interior Neumann problem for Laplace's equation the kernel $k_1$ is smooth whereas $k_2$ is a combination of a smooth function and a logarithmic kernel, i.e.
\begin{equation}\label{kernels}
k_2(x,y)=\rho(x,y)+\psi(x) \log|x-y|,
\end{equation}
with $\rho$ and $\psi$ smooth functions.

Motivated by these applications, in this paper we treat \eqref{EqInFr}
when the kernels are continuous functions and/or weakly singular at the bisector as, for example, $|x-y|^\nu$, $\nu>-1$ and $\log|x-y|$. Without losing the generality, we first consider the case when $k_1$ and $k_2$ are both smooth and the case when $k_1$ is smooth and $k_2$ is weakly singular. However, our approach can be also applied in other ``mixed'' situations, as, for instance, \eqref{kernels}.

Let us note that if $k_1(x,y) \equiv 0$, then equation \eqref{EqInFr} is the classical nonlinear Hammerstein equation
\begin{equation*}
f(y) -  \int_{-1}^{1} k_2(x,y)h(x,f(x)) dx= g(y), \qquad y \in [-1,1],
\end{equation*}
which is one of the most frequently investigated nonlinear integral equations, since it occurs in applications in  numerous areas. Several problems written in terms of ordinary and partial differential equations can be transformed into equations of Hammerstein type through  Green's function  \cite{Atkinsonsurvey,Pascali,Wazwaz}. An example is the following differential problem of the second order which describes the forced oscillations of finite amplitude of a pendulum  \cite{Chidume}
\begin{equation}\label{difproblem}
\begin{cases}
F_{yy}(y)+a^2 \sin{F(y)}=G(y), \quad y \in [0,1] \\
F(0)=F(1)=0,
\end{cases}
\end{equation}
where $F$ denotes the amplitude of oscillation, the constant $a \neq 0$ depends on the length of the pendulum and on gravity, and the driving force $G$ is periodic and odd. Problem \eqref{difproblem} is equivalent to this nonlinear integral equation
$$F(y)+\int_0^1 k(x,y) \left[G(x)-a^2 \sin{(F(x))}\right] dx=0, $$
with $k(x,y)$ the triangular function defined as
$$k(x,y)= \begin{cases}
x(1-y), & x \in [0,y] \\
y(1-x), & x \in [y,1]
\end{cases}.$$
A further example is the well-known Chandrasekhar H-equation
$$H(y)-c H(y) \int_{0}^1 \frac{x s(x) }{y+x}  H(x) dx=1, \qquad c \in \mathbb{C},$$
where $s$ is a given function and $H$ is the unknown. It models various physical problems such as the radioactive transfer and the  kinetic of gases and, setting $f(y)=[H(y)]^{-1}$, it can be written as an Hammerstein equation \cite{Atkinsonsurvey} in the unknown $f(y)$. \newline
Other contexts of application are the network theory, optimal control systems and automation \cite{Dolezale,Narendra}.

Many interesting papers on the approximation of the solution of Hammerstein equations have appeared in the last few years. The survey \cite{Atkinsonsurvey} provides a complete overview of methods that can be also applied to other kind of nonlinear integral equations. Detailed examples of existing methods are collocation methods \cite{Kumar90,KumarSloan}, degenerate kernel methods \cite{Kaneko91}, discrete Legendre spectral methods \cite{das2016} also for weakly singular equations \cite{Mandal21,Mandal1}, and the more recent  numerical techniques based on spline quasi-interpolation \cite{barrera20,dagnino19} and Gaussian spline rules \cite{barrera22}.

In this paper, first we determine the functional spaces which the solution of the equation belongs to, and study the mapping properties of the involved integral operators by using suitable approximation tools. Then, we propose Nystr\"om type methods based on the polynomial approximation. Although this approach has been widely applied to linear Fredholm integral equations of the second kind (see, for instance, \cite{debonis06,FermoRusso2010,FermoLaurita2015,Laurita2020}), this is the first time that it is developed for nonlinear second-kind equations, according to our knowledge.
The Nystr\"om method is based on a discretization of the integral operators which involves the Gauss-Legendre rule if the kernel is smooth or a suitable product rule, based on the Legendre nodes, if the kernel is weakly singular. Following \cite{KEAtkinson}, we prove the stability of the method in spaces equipped with the uniform norm and we provide new estimates of the error, deduced also  thanks to the recent results given in \cite{BestApproxErr}. Specifically, under suitable assumption on the known functions, we prove that the rate of convergence of the method is comparable with the error of best polynomial approximation in the functional spaces where the solution lives.  

We conclude the paper by applying the proposed method to the numerical solution of the interior Neumann problem for the Laplace equation having nonlinear boundary conditions, following an approach which has been already shown in \cite{AtkChan,ruotsalainen1988}.

The paper is structured as follows. Section \ref{sec:preliminaries}
details the function spaces in which equation \eqref{EqInFr} is considered and provides some basic results concerning the error of best polynomial approximation. Section \ref{sec:operators} focuses on the mapping properties of the involved integral operators and on the solvability of equation \eqref{EqInFr}. Section \ref{sec:smooth} and Section \ref{sec:weakly} deal with the Nystr\"om methods we propose when the kernels are smooth or weakly singular, respectively. In both situations, a theoretical study is provided together with some numerical experiments which show  the performance of the method. Section \ref{sec:BIE} concerns an application of the described numerical approach to the Laplace equation with nonlinear Neumann boundary conditions. Section \ref{sec:proofs} contains the proofs of our results.

\section{Function spaces and best polynomial approximation}\label{sec:preliminaries}
Let us denote by $C^0 \equiv C^0([-1,1])$ the Banach space of continuous functions on $[-1,1]$ with the uniform norm
\[ \|f\|_{C^0} = \|f\|_\infty= \max_{x \in [-1,1]}|f(x)|,\]
and let us  introduce the  Sobolev--type space of order $1 \leq r \in \mathbb{N}$
\begin{equation*}
{W}^r=\left\{f \in C^0: f^{(r-1)} \in AC((-1,1)), \,  \|f^{(r)}\varphi^r  \|_\infty <\infty \right\},
\end{equation*}
where $\varphi(x)=\sqrt{1-x^2}$ and $AC((-1,1))$ denotes the set of all absolutely continuous functions on $(-1,1)$.
We equip $W^r$ with the norm
$$\|f\|_{W^r}=\|f\|_\infty+\|f^{(r)}\varphi^r \|_\infty.$$
The space $(W^r, \| \cdot\|_{W^r})$ is a Banach space.

For our aims, we also need to define the Sobolev space ${\boldsymbol W}^{r}(\Omega)$ for bivariate functions $f:\Omega \rightarrow \mathbb{R}$, with $\Omega$ open subset of $ \mathbb{R}^2$. It is
the set of all functions $f$ in $\Omega$ such that for every $2$-tuple of nonnegative integers $\ell=(\ell_1,\ell_2)$, with $\displaystyle |\ell|=\sum_{i=1}^2 \ell_i \leq r$, the mixed partial derivatives
$D^{\ell} f = \frac{\partial^{\ell_1+\ell_2} f}{\partial x_1^{\ell_1} \partial x_2^{\ell_2}}$
exist and $\|D^{\ell} f\|_\infty <\infty$ .
We endow this space with the norm
\begin{equation*}
||f||_{{\boldsymbol W}^{r}(\Omega)}=||f||_\infty+ \sum_{1 \leq |\ell| \leq r} \|D^{\ell} f\|_\infty.
\end{equation*}
For functions of ``intermediate'' smoothness, we define the Zygmund space $Z^\la$,  with  $\la \in \RR^+$, as follows
\[Z^\la=\left\{f\in C^0 \ : \ \sup_{t> 0}\frac{\Omega_{\vphi}^k(f,t)}{t^\la} <\infty, \quad  k \geq 1, \,  k> \la \right\},\]
where the main part of the \textit{$\varphi$-modulus of smoothness} $\Omega_{\vphi}^k(f,t)$ is defined as  \cite[p. 90]{MastroianniMilovanovic}
\begin{equation}\label{Omega}
\Omega_{\vphi}^k(f,t)=\sup_{0<\tau\leq t}\max_{x\in I_{k\tau}}|\Delta_{\tau\vphi}^k f(x)|,\quad I_{k\tau}=[-1+(2k\tau)^2,1-(2k\tau)^2],
\end{equation}
with
\[\Delta_{\tau\vphi}^k f(x)=\sum_{i=0}^k (-1)^i \binom k i f\left(x+\frac{\tau\vphi(x)}2(k-2i)\right).\]
The space $Z^\la$ is endowed with the  norm
\begin{eqnarray*}
\|f\|_{Z^\la}&=&\|f\|_\infty +  \sup_{t> 0}\frac{\Omega_{\vphi}^k(f,t)}{t^\la},
\end{eqnarray*}
and also ($Z^\la, \|\cdot\|_{Z_\la}$) is a Banach space.

From now on we will denote by $\mathcal{C}$ a generic positive constant that can be different in different formulas. Moreover will we write $\C =\C(a,b,...)$ to say that
 $\C$ is dependent on the parameters $a,b,....$ and $\C \neq \C(a,b,...)$ to say that $\C$ is independent of them.

Denoting by $\PP_m$ the set of all algebraic polynomials of degree at most $m$,
for functions $f \in C^0$, let us now define  the error of best polynomial approximation as
\[E_m(f)=\inf_{P_m \in \PP_m} \left\|\left(f-P_m\right)\right\|_\infty.\]

It is well known that (see, for instance, \cite{MastroianniMilovanovic} and the references therein)
\[
f\in C^0 \ \Longleftrightarrow \ \lim_{m\rightarrow \infty} E_m(f)=0,
\]
Moreover it is known that the behavior of the best approximation error is strictly related to the smoothness of the function $f$.
Indeed in order to estimate $E_m(f)$, we can use, for instance, the following weak-Jackson inequality \cite{DT}
\begin{equation}
	\label{jackson}
	E_m(f)\le \C \int_0^{\frac 1 m}\frac{\Omega_{\vphi}^k(f,t)}t dt,\quad \forall f\in C^0, \quad \C\neq \C(m,f).
\end{equation}
On the other hand, it is also well known the so called Favard inequality, which says that
\begin{equation}\label{erroreSobolev}
	E_m(f) \leq \frac{\mathcal{C}}{m^r}\|f\|_{W^r}, \qquad \forall f \in W^r,
\end{equation}
where $\mathcal{C}\neq \mathcal{C}(m,f)$.

A stronger result is that both the Sobolev and Zygmund spaces introduced before, can be characterized in terms of the best polynomial approximation error. Indeed, from \cite[Th. 4.2.1, p. 40]{DT} and \cite[Co.2.2, p. 224]{LR} it immediately follows that
\begin{equation}\label{Ch.Sobolev}
f\in W^r \Longleftrightarrow  E_m(f) =\mathcal{O} \left(\frac{1}{m^r}\right), \qquad 	\forall r\in \mathbb{N}, 
\end{equation}
while by \cite[Th.8.2.1., p. 94]{DT} and (\ref{jackson}) it follows
\begin{equation}\label{Ch.Zyg}
 f\in Z^\lambda \Longleftrightarrow  E_m(f) =\mathcal{O} \left(\frac{1}{m^\lambda}\right), \qquad 	\forall \lambda \in \mathbb{R^+}\setminus\mathbb{N}.
\end{equation}

We conclude the section with a very recent result \cite{BestApproxErr} which provides a characterization of the error of the best polynomial approximation of composite functions and that will be useful in the sequel.
\begin{theorem} \label{teo:BestApprx}
Let $h: \Omega \rightarrow \mathbb{R},$ with $\Omega$ open subset of $\mathbb{R}^2$ and $\sigma: (-1,1) \rightarrow \mathbb{R}^2$ such that $Im( \sigma) \subseteq \Omega.$
Assume that $h \in \mathbf{W}^r( \Omega)$ and $ \sigma(x)= (x, f(x))$ with $ f \in W^r,$
then
\[
E_m( h \circ \sigma) \leq  \mathcal{C} \left(\frac{2}{m}\right)^r B_r
\, \|h\|_{\mathbf{W}^r( \Omega)} \, \|f\|_{W^r}^s,
\]
where $\mathcal{C}=\mathcal{C}(r)$ is a positive constant independent of $h$ and $f$, $B_r$ is the $r$-th Bell number,  and the exponents $s$ is defined as follows
\begin{equation}\label{esse}
s=\left\{\begin{array}{lcr} 0, & \quad \mathrm{if}  \, & \|f\|_{W^r} \leq 1 \vspace{0.1cm} \\
r, & \quad  \mathrm{if}  \,& \|f\|_{W^r} >1  \end{array} \right..
\end{equation}
\end{theorem}
 
 \section{The solvability of equation (\ref{EqInFr})}\label{sec:operators}
The aim of this section is, firstly, to investigate the mapping properties of the operators involved in equation \eqref{EqInFr} and then to study its solvability in suitable subspaces of $C^0$.

We start with the case that both $k_1(x,y)$ and $k_2(x,y)$ are smooth functions.

Let us introduce the Fredholm operators
\begin{equation}\label{K}
(K^if)(y)= \int_{-1}^{1} k_i(x,y) f(x) dx, \qquad y \in [-1,1], \qquad i=1,2,
\end{equation}
and the so-called Nemytskii operator
\begin{equation}\label{H}
(Hf)(x)= h(x,f(x)), \qquad x \in [-1,1].
\end{equation}
Then, equation \eqref{EqInFr} can be written as
\begin{equation}\label{Hamm_op}
(I-\mathcal{K})f=g, \qquad \mathcal{K}= K^1+ K^2H
\end{equation}
where $I$ is the identity operator.

It is well-known that if we assume the following hypothesis
\begin{itemize}
\item[\bf{[K1]}] The kernels $k_i(x,y)$, $i=1,2$, are such that
\begin{align*}
&\sup_{y \in [-1,1]} \int_{-1}^1 |k_i(x,y)| dx<\infty, \\ &\lim_{y \to \tilde{y}} \int_{-1}^1 |k_i(x,y)-k_i(x,\tilde{y})| dx=0, \quad \tilde{y} \in [-1,1];
\end{align*}
\end{itemize}
then the linear Fredholm operators $K^i:C^0 \to C^0$, $i=1,2$, are compact, from which we can also deduce that they are completely continuous; see, for instance, \cite{A}.

Moreover, if we assume that
\begin{itemize}
\item[\bf{[H1]}] The function $h:[-1,1] \times \mathbb{R} \to \mathbb{R} $ is continuous;
\item[\bf{[H2]}] The partial derivative $h_v(x,v)= \dfrac{\partial h(x,v)}{\partial v}$
exists and is continuous for $x \in [-1,1]$ and $v \in \mathbb{R}$;
\end{itemize}
then the Nemytskii operator $H:C^0 \to C^0$ is well defined, bounded, and continuous because of the hypothesis {\bf{[H1]}} (see, for instance, \cite{Pot2022}) and is continuously Fr\'echet differentiable on the space $C^0$ thanks to {\bf{[H2]}} (see, for instance, \cite[Lemma 4]{KumarSloan}).

Its Fr\'echet derivative, at $f\in C^0$, is given by the multiplicative linear operator defined as
$$[(H'f) \phi](x)= h_v(x,f(x)) \phi(x),  \qquad \forall \phi \in C^0, \quad x\in[-1,1]. $$

For our aims, it is also useful introduce the linear operator
\begin{equation*}
(Gf)(y)=(K^2 f)(y)+ g(y), \qquad  y \in [-1,1].
\end{equation*}
Let us note that if the further condition
\begin{itemize}
\item[\bf{[G1]}] The right-hand side $g$ is continuous in $[-1,1]$;
\end{itemize}
is fulfilled, the operator $G$ inherits the same properties as $K^2$.
Consequently, the composite operator $GH:C^0 \to C^0 $ is completely continuous. \newline

Let us now rewrite equation \eqref{EqInFr} (or equivalently \eqref{Hamm_op}) as the following fixed point problem
\begin{equation}\label{EqInFrOp}
 f(y)= (\mathcal{G}f)(y), \qquad (\mathcal{G}f)(y)= (K^1f)(y)+(GH f)(y)=(\mathcal{K}f)(y)+g(y).
\end{equation}

Under the assumptions {\bf{[K1]}}, {\bf{[H1]}}, {\bf{[H2]}}, and {\bf{[G1]}}
the operator $\mathcal{G}$ is continuously Fr\'echet differentiable on the space $C^0$. For each $f \in C^0$, its derivative is given by
\begin{align*}
[(\mathcal{G}'f)\phi](y) &=[( K^1+(GH)'f)\phi](y) \\ & =  (K^1\phi)(y)+ [G(H'f)\phi](y), \qquad \forall \phi \in C^0, \quad y\in [-1,1],
\end{align*}
that is
\begin{align}\label{G'}
[(\mathcal{G}'f)\phi](y)  =  \int_{-1}^1 k_1(x,y) \phi(x) dx +\int_{-1}^{1}{k_2(x,y)h_{v}(x,f(x))\phi(x)}dx.
\end{align}

Moreover, if
\begin{itemize}
\item[\bf{[H3]}] The partial derivative $h_{vv}(x,v)= \dfrac{\partial^2 h(x,v)}{\partial^2 v}$
exists and is continuous for $x \in [-1,1]$ and $v \in \mathbb{R}$;
\end{itemize}
the operator $\mathcal{G}$ also admits  the second  Fr\'echet derivative at $f\in C^0$, defined as \cite{KEAtkinson}
\begin{equation}\label{G''}
[(\mathcal{G}''f)(\phi_1,\phi_2)](y)=  \int_{-1}^{1}{k_2(x,y)h_{vv}(x,f(x))\phi_1(x)\phi_2(x)}dx,   \quad \phi_1,\phi_2 \in C^0.
\end{equation}

\begin{remark}
Let us observe that if $h \in \mathbf{W}^r(\Omega)$ where $\Omega$ is an open subset of $\mathbb{R}^2$ and $r \geq 2$, then the assumptions {\bf{[H1]}}, {\bf{[H2]}}, and {\bf{[H3]}}  are satisfied.
\end{remark}

The analysis of the operator $\mathcal{G}$ is fundamental to state the existence of solutions of (\ref{EqInFrOp}).
Let us recall that a solution $f^*$ of (\ref{EqInFrOp}) is \textit{ geometrically isolated } if there exists a ball
\begin{equation*}
B(f^*, \delta)= \lbrace f \in C^0 : \|f-f^*\|_\infty \leq \delta \rbrace,
\end{equation*}
for some $\delta>0$ that does not contain any solution of \eqref{EqInFrOp} other than $f^*$. 

Let us also remind that the \textit{index} of a geometrically isolated solution $f^*$ is the common value of the \textit{rotation} of the vector field $I- K^1+GH$ over all sufficiently small spheres centered at $f^*$; see, for instance, \cite[p. 100]{Krasnoseski} and \cite[Section 2]{KumarSloan}.

The existence and the uniqueness of a geometrically isolated and with nonzero index solution $f^*$ of \eqref{EqInFrOp} is established in the next theorem, which can be deduced by the more general result  \cite[Theorem 21.6, p. 108]{Krasnoseski}.

\begin{theorem}\label{TeoEsUn}
Assume that the operator $ \mathcal{G}: C^0\rightarrow  C^0$ defined in \eqref{EqInFrOp} is  completely continuous. Let $f^*$ be such that  $I-\mathcal{G}'(f^*)$ is invertible,  where $\mathcal{G}'(f^*)$ is the Fr\'{e}chet derivative at the point $f^*$.
Then $f^*$ is a fixed point of $\mathcal{G}$.  Moreover, assume that $1$ is not an eigenvalue of $\mathcal{G}'(f^*).$  Then $f^*$ is the unique nonzero index geometrically isolated solution of equation \eqref{EqInFrOp} in $C^0$.
\end{theorem}

About the smoothness of the solution we can state the following result.

\begin{theorem}\label{teo:K}
	Let $\mathcal{K}=K^1+K^2H$ where $K^1$, $K^2$, and $H$ are defined in \eqref{K} and \eqref{H}, respectively.  Under the assumptions of Theorem \ref{TeoEsUn} if
	$$\sup_{x \in [-1,1]} \|k_i(x,\cdot)\|_{W^r}< \infty, \qquad i=1,2, $$
	then $\forall f\in C^0$ it is $\mathcal{K}f\in W^r$.
	Consequently, if also $g \in W^r$ the solution $f^*$ of \eqref{EqInFrOp} belongs to $W^r$.
\end{theorem}
Now, we consider the case when  equation \eqref{EqInFr} presents a nonlinear operator having a symmetric weakly singular kernel at the bisector, i.e.
\begin{equation}\label{EqInFrW}
	f(y) -  \int_{-1}^{1} k_1(x,y)f(x) dx-  \int_{-1}^{1} k_2(x,y) h(x,f(x)) dx= g(y), \qquad y \in [-1,1],
\end{equation}
where $$k_2(x,y)=\psi(x) k^*(|x-y|)$$ with $\psi$ a smooth function on $[-1,1]^2$. Typical examples of the kernel $k^*(|x-y|)$  are
\begin{equation}\label{k*}
	k^*(|x-y|)= |x-y|^\mu, \quad \mu>-1, \qquad k^*(|x-y|)= \log{|x-y|}.
\end{equation}
About the solvability of this equation we have the following result.
\begin{theorem}\label{theoWeak}
	Assume that   {\bf{[H1]}}, {\bf{[H2]}}, and {\bf{[G1]}} holds true. If
\begin{itemize}
	\item[\bf{[H4]}]
\begin{align*}
\hspace*{-1cm}\max_{y\in [-1,1]} & \left(\int_{-1}^{1}|k_1(x,y)|dx+\max_{x\in[-1-1]}\|h_v(x,\cdot)\|_{\infty}  \int_{-1}^{1}|\psi(x)k^*(|x-y|)|dx  \right)<1
\end{align*}	
	\end{itemize}
	then equation (\ref{EqInFrW}) has a unique solution $f^*\in C^0$.
	\end{theorem}
	
Assumption   {\bf{[H4]}} is quite restrictive. Nevertheless in the case (\ref{k*}) it is possible to show that  {\bf{[K1]}} is still satisfied (see \cite{A}) and therefore Theorem \ref{TeoEsUn} is still true for this special, but frequent, case.  

Moreover in this special case, we can also state a result about the smoothness of the solution of the equation \eqref{EqInFrW}. 

In  what follows we adopt the notation
	\begin{equation}\label{Kmu}
		(\widetilde{K}^\mu f)(y)=\left\{\begin{array}{lc}\displaystyle \int_{-1}^1 \psi(x) |x-y|^\mu f(x) dx, & \mu>-1, \quad \mu\neq 0, \\

\displaystyle 	\int_{-1}^1 \psi(x) \log|x-y| f(x) dx, & \mu=0	\end{array}\right..
	\end{equation}
	
	\begin{theorem}\label{classeDebSing}
	Let $\mathcal{K}=K^1+\widetilde{K}^\mu H$, $\mu>-1$, where $K^1$, $\widetilde{K}^\mu$, and $H$ are defined in \eqref{K}, \eqref{Kmu} and \eqref{H}, respectively.  	Assume that equation (\ref{EqInFrW}) has a unique solution $f^*\in C^0$. If  {\bf{[H1]}} holds, $\psi\in C^0$, and  
\begin{equation}\label{assk1}
\sup_{x\in[-1,1]}\|k_1(x,\cdot)\|_{Z^\lambda}<\infty, \quad \lambda >0,
\end{equation}	
then $\forall f\in C^0$ it is $\mathcal{K}f\in Z^{s}$, where $s=\min\{\lambda, \mu+1\}$, when $\mu\neq 0$, while $s=\min\{\lambda, 1-\epsilon\}$, with $\epsilon>0$ arbitrarily small, when $\mu = 0$.  Consequently, if also $g \in Z^{s}$ the solution $f^*$ of the equation $(I-\mathcal{K})f=g$ belongs to $Z^{s}$.
\end{theorem}
\begin{remark}
In the literature some estimates of the smoothness of the solution of Hammerstein integral equations with a weakly singular integral of the type (\ref{Kmu}) are known (see, for instance, \cite{Kaneko90} and  \cite{PV}). Here, for a more general equation, we give minimal smoothness assumptions on the known functions, in order to determine the space which the solution of the equation belongs to. 
\end{remark}

\section{The Nystr\"om method for the case of continuous kernels}\label{sec:smooth}
We introduce now a numerical method of Nystr\"om type based on the Gauss-Legendre rule
\begin{equation}\label{gauss}
\int_{-1}^1 f(x) dx = \sum_{k=1}^m \lambda_k f(x_k)+e_m(f),
\end{equation}
where $\lambda_k$ is the $k$th Christoffel number, $x_k$ is the $k$th zero of the orthonormal Legendre polynomial $p_m(x)$ of degree $m$, and $e_m(f)$ is the remainder term.
Let us recall that \cite[Theorem 5.1.6]{MastroianniMilovanovic}
\begin{equation}\label{em}
|e_m(f)| \leq 4 \, E_{2m-1}(f), \qquad \forall f \in C^0.
\end{equation}

First, let us approximate the Fredholm integral operators (\ref{K}) by defining the discrete operator $K^i_m$ as follows
\begin{equation}\label{Km}
(K^{i}_{m}f)(y)= \sum_{k=1}^{m} \lambda_{k} k_i(x_k,y) f(x_k), \qquad i=1,2.
\end{equation}
Now, consider the equation
\begin{equation}\label{eqFinDim}
(I-\mathcal{K}_m)f_m=g, \qquad \mathcal{K}_m= K^1_m+K^2_mH
\end{equation}
i.e.
\begin{equation}\label{Gmfm}
f_m= \mathcal{G}_m f_m, \qquad \text{with} \quad \mathcal{G}_mf=\mathcal{K}_mf+g
\end{equation}
where $f_m$ is an unknown function, approximating the solution of equation \eqref{EqInFrOp}.

The next theorem contains some useful properties of the operators $\mathcal{K}_m$  essential for obtaining the stability and convergence of the method. Here and in the following all the involved operators will be considered as maps of $C^0$ into itself.
\begin{theorem}\label{teo:K_Km}
Let $\mathcal{K}$ and $\mathcal{K}_m$ be the operators defined in \eqref{Hamm_op} and \eqref{eqFinDim}, respectively, and let us assume that ${\bf{[K1]}}$ and ${\bf{[H1]}}$ are satisfied.  Then, the sequence $\{\mathcal{K}_m\}_m$ is collectively compact  and pointwise convergent to $ \mathcal{K}$.
\end{theorem}

\begin{remark}
Note that if, in addition to ${\bf{[K1]}}$ and ${\bf{[H1]}}$, the  assumption ${\bf{[G1]}}$ is satisfied, from Theorem \ref{teo:K_Km} we can deduce that the sequence $\{\mathcal{G}_m\}_m$ is collectively compact  and pointwise convergent to $ \mathcal{G}$.
\end{remark}

Now, in order to compute the unknown solution $f_m$ of equation \eqref{eqFinDim}, which has the explicit form
\[ f_m(y)-  \sum_{k=1}^{m} \lambda_k k_1(x_k,y) f_m(x_k)- \sum_{k=1}^{m}{ \lambda_k k_2(x_k,y) h(x_k, f_m(x_k))}=g(y),\]
let us collocate it at the points $x_i$ for $i=1,\dots,m$. In this way, we obtain the following nonlinear system of $m$ equations in the $m$ unknowns $a_i= f_m(x_i),$ $i=1,\ldots,m$,
\begin{align}\label{NonLinSis}
\sum_{k=1}^{m}[\delta_{i,k} -  \lambda_{k} k_1(x_k, x_i)] a_k-\sum_{k=1}^{m}{ \lambda_{k} k_2(x_k, x_i) h(x_k, a_k)}&= g(x_i),\\ &  i=1,\ldots,m, \nonumber
\end{align}
where $\delta_{i,k}$ is the Kronecker symbol.

The solution $(a^*_1,\ldots,a^*_m)$ of system (\ref{NonLinSis}), allows us to construct the\\
Nystr\"om interpolant
\begin{equation}\label{interpolante}
f_m(y)=\sum_{k=1}^{m}{\lambda_k [k_1(x_k,y) a^*_k+k_2(x_k,y) h(x_k,a^*_k)]}+g(y).
\end{equation}

\begin{theorem} \label{teo:stability}
Assume {\bf{[K1]}}, {\bf{[H1]}}, {\bf{[G1]}}, and {\bf{[H2]}} and {\bf{[H3]}} in the ball $B(f^*,\delta)$, where $f^*$ is a fixed point of the operator $\mathcal{G}$ defined in \eqref{EqInFrOp}. Moreover, assume that  $1$ is not an eigenvalue of $\mathcal{G'}f^*$ where $\mathcal{G'}f^*$ is given in \eqref{G'}.
Then, for $m$ sufficiently large, say $m \geq m_0$, the operator $\mathcal{G}_m$ in \eqref{Gmfm} has a unique fixed point $f_m$ in $B(f^*, \epsilon)$ with $0< \epsilon\leq \delta$.
\end{theorem}

\begin{theorem}\label{teo:conv}
Assume that the hypotheses of Theorem \ref{teo:stability} are fulfilled.
Let $f^*$ be the unique fixed point of $\mathcal{G}$ in $B(f^*,\delta)$, for some $\delta>0$. Assume that $1$ is not an eigenvalue of $\mathcal{G'}(f^*)$ and, in addition, that for some $r\geq 1$,
\begin{equation}\label{assK}
 \sup_{x \in [-1,1]} \|k_i(x,\cdot)\|_{W^r}< \infty,  \qquad \sup_{y \in [-1,1]} \|k_i(\cdot,y)\|_{W^r}< \infty,\qquad i=1,2
\end{equation}
\begin{equation}\label{assgh}
g \in W^r, \qquad h \in \mathbf{W}^r( \Omega), 
\end{equation}
where $ \Omega$ is an open subset of $\mathbb{R}^2$.

Then, for $m$ sufficiently large (say $m \geq m_0$), denoted by $f_m$ the unique fixed point of $\mathcal{G}_m$ in $B(f^*,\epsilon)$, with $0 < \epsilon \leq \delta$, we have
\begin{equation}\label{StimaErroref}
||f^* -f_m||_\infty = \mathcal{O} \left(  \frac{1}{m^r} \right).
\end{equation} 
\end{theorem}
\subsection{Numerical tests}
In this section, we consider some numerical tests to confirm the effectiveness of our method.

In all the experiments, first we solve system \eqref{NonLinSis} by using the classical Newton method or the Matlab routine \texttt{fsolve}. Then, we compute the Nystr\"om interpolant $f^*_m$ given in \eqref{interpolante}, and the relative discrete errors on a grid of $10^2$ equidistant nodes $y_i$, $i=1,\dots,100$, in $[-1,1]$, i.e.
\begin{equation}
\mathcal{E}_m= \frac{\|f^{*} - f^{*}_{m}\|}{ \|f^{*}\|}, \qquad
\label{RelDErr}
\end{equation}
where $\|f\|=\displaystyle  \max_{i=1,\dots,10^2} |f(y_i)|$ and  $f^*$ is the exact solution.

The first two examples, in which the exact solution is known, aim to make a comparison with other methods available in the literature.
In the other two, the solution $f^{*}$ is not known, and then we consider as exact the Nystr\"om interpolant $f^*_{512}$. Moreover, in these specific tests, the known functions are not so smooth so that we can show the efficacy of the method also in these cases.

In addition,  we compute the estimated order of convergence 
\begin{equation}\label{EOC}
EOC_m=\frac{\log{(\mathcal{E}_m/\mathcal{E}_{2m})}}{\log 2}.
\end{equation}

\begin{example}\label{Example1}
\rm
Consider the Hammerstein equation
\[ f(y)- \int_{-1}^{1}{ e^{y-2x} (f(x))^3}dx= e^{y-1}( e- e^{2} +1), \quad \vert y\vert \leq 1,
\]
whose exact solution is $f(x)= e^{x}$. Such equation has been considered in \cite{barrera22,Shav}. In \cite{barrera22} the machine precision is achieved by solving a nonlinear system of order $64$,  whereas in \cite{Shav} the better convergence order is $10^{-10}$.     Table \ref{table:Example1} shows that our method has a faster convergence. In fact, the machine precision is reached by solving a nonlinear system of $m=8$ equations in $6$ iterations. This is certainly due to the smoothness of the known functions which are analytic in $[-1,1]$.

\begin{table}
\caption{From left to right the numerical results for Example \ref{Example1}, and Example \ref{Example2} \label{table:Example1}}
\begin{tabular}{lll}
\hline\noalign{\smallskip}
				$m$ &  $\mathcal{E}_m$ & iter\\
\noalign{\smallskip}\hline\noalign{\smallskip}
				4 	& 4.88e-08 &  6\\
				8 	& 4.90e-16 & 6\\
\noalign{\smallskip}\hline
\end{tabular}
\qquad 
\begin{tabular}{lll}
\hline\noalign{\smallskip}
				$m$ &  $\mathcal{E}_m$ & iter\\
\noalign{\smallskip}\hline\noalign{\smallskip}
				4 		& 4.87e-03 & 5  \\
				8 	& 2.32e-07 & 5	 \\
				16 		 & 2.22e-16 & 5 \\
\noalign{\smallskip}\hline
\end{tabular}
\end{table}
\end{example}

\begin{example}\label{Example2}
\rm
Let us  apply our method to the  Hammerstein equation \cite{RKP}
\[ f(y)-\int_{-1}^{1}{ y \cos{ \left(\frac{\pi}{2}x \right)}e^{f(x)}}dx= \sin{ \left(\frac{\pi}{2}y \right)} - \frac{4y}{\pi} \sinh{(1)}\quad \vert y\vert \leq 1,
\]
where the right-hand side term is fixed so that the exact solution is $f(x)= \sin{ \left(\frac{\pi}{2}y \right)}$. Also in this case, the kernel and the right-hand side are analytic functions and then we expect a fast convergence. Looking at the errors (see Table \ref{table:Example1}) we can note that  our method reaches the machine precision with $m=16$. We remark that the best convergence error of the method presented in \cite{RKP}  is equal to $10^{-8}$.
\end{example}

\begin{example}\label{Example3}
\rm In this test we consider the complete equation of the form \eqref{EqInFr}
\[ f(y)-\int_{-1}^{1}{y \cos{x}}dx-\int_{-1}^{1}{\frac{e^{x+y} \cos{(x+1)}}{x^2+5} \frac{dx}{1+(f(x))^2}}= |y|^{\frac{5}{2}}, \quad \vert y\vert \leq 1.
\]
where kernels are both smooth and  the right hand side term $g \in W_2$. Then, according to \eqref{StimaErroref}, the expected theoretical order of convergence is $\mathcal{O}(m^{-2})$. The numerical results reported in Table \ref{table:Example3} (on the left) shows a better convergence.
  
\begin{table}
\caption{From left to right the numerical results for Example \ref{Example3}, and  Example \ref{Example4}  \label{table:Example3}}\begin{tabular}{llll}
\hline\noalign{\smallskip}
				$m$ &  $\mathcal{E}_m$ & iter & $EOC_m$\\
\noalign{\smallskip}\hline\noalign{\smallskip}
				8 	& 2.35e-04 & 4	& 3.85e+00  \\
                16 	& 2.15e-05 & 4 & 3.45e+00	\\
                32 	& 1.98e-06 & 4 &  3.44e+00  \\
               64 	& 1.79e-07 & 4 	&  3.47e+00 \\
              128 	& 1.59e-08 & 4 	& 3.49e+00 \\
               256 	& 1.30e-09 & 4 	& 3.61e+00 	  \\
\noalign{\smallskip}\hline
\end{tabular}
\qquad 
\begin{tabular}{llll}
\hline\noalign{\smallskip}
				$m$ &  $\mathcal{E}_m$ & iter & $EOC_m$\\
\noalign{\smallskip}\hline\noalign{\smallskip}
				 8 & 9.45e-04  & 18 & 7.96e+00 	 \\
               16 &	 4.77e-05 	& 65 &	4.31e+00  \\
               32 	& 2.26e-06 	& 21 &	4.40e+00  \\
               64 	& 1.03e-07 	& 20 &  4.45e+00 \\
              128 &	 4.64e-09 	& 20 & 4.48e+00 \\
              256 &	 1.98e-10 & 20 	& 4.55e+00  \\
\noalign{\smallskip}\hline
\end{tabular}
\end{table}

\end{example}
\begin{example}\label{Example4}
\rm
Consider an equation in which the kernel of the nonlinear operator satysfies conditions \eqref{assK} with $r=3$
\[ f(y)-\int_{-1}^{1}{(x+y)f(x)}dx-\int_{-1}^{1}{|xy|^{\frac{7}{2}} (f(x))^3}dx= e^y +\log{(3+y)}, \quad \vert y\vert \leq 1.
\]
Table \ref{table:Example3} (on the right) shows, also in this case, that the numerical errors  are better than the expected once, since the theoretical error is of the order $\mathcal{O}(m^{-3})$.
\end{example}
	
\section{The Nystr\"om method for the case of weakly singular kernels} \label{sec:weakly}
Let us introduce a product rule which allows us to approximate integrals having kernels with weak singularity at the bisector of the following type
\[\mathcal{I}(f,y)=\int_{-1}^{1}\psi(x) k^*(|x-y|) f(x)dx.\]
The quadrature formula is given by
\begin{equation}\label{productrule}
\mathcal{I}(f,y)= \sum_{k=1}^{m}{\left[\lambda_k \sum_{i=0}^{m-1} p_i(x_k) M_i(y)\right]f(x_k)}+r_{m}(f,y)= \mathcal{I}_m(f,y)+r_{m}(f,y),
\end{equation}
where $\lambda_k$ and $x_k$ are the $k$th Christoffel number and the $k$th zero of the Legendre polynomial $p_m$, respectively, $M_i(y)$, $i=0,1,\ldots,m-1$, are the so-called modified moments defined as
\begin{equation}\label{modmom} M_i(y)=\int_{-1}^{1}{p_i(x)\,  \psi(x) k^*(|x-y|)} \, dx,  \qquad i=0,\dots,m-1,
\end{equation}
and $r_m(f,y)$ is the quadrature error.
From now on, we will denote the weights of rule \eqref{productrule} by
$$c_k(y)=\lambda_k \sum_{i=0}^{m-1} p_i(x_k) M_i(y), \qquad k=1,\ldots,m.$$ 

Next theorem provides the assumptions assuring the stability of the \\ quadrature formula \eqref{productrule} and an estimate for the remainder term (see \cite[Theorem 5.1.11]{MastroianniMilovanovic} and \cite{Mezzanotte2021}).
\begin{theorem}
Let us assume that the kernel $k^*$ and the function $\psi$ satisfies
\begin{align*}
& \sup_{|y| \leq 1} \int_{-1}^1 \frac{|\psi(x)k^*(|x-y|) |}{\sqrt{\varphi(x)}} dx < \infty, \\
& \sup_{|y| \leq 1} \int_{-1}^1 \psi(x)k^*(|x-y|)  \left(1+ \log^{+} \psi(x)k^*(|x-y|)\right) \, dx \, <\infty,
\end{align*}
where $\log^{+} \psi(x) k^*(|x-y|)=\log{\max\{1,\psi(x) k^*(|x-y|) \}}$. Then,
for each $f \in C^0$ one has
\[  \sup_m{\sup_{ |y| \leq 1}{|\mathcal{I}_m(f,y)|}} \leq \mathcal{C} \|f\|_\infty,
\]
and
\begin{equation}\label{prodotto}\sup_{ |y| \leq 1}{|r_{m}(f,y)|} \leq \mathcal{C} E_{m-1}(f),
\end{equation}
with $\mathcal{C}\neq \mathcal{C}(m,f).$
\end{theorem}

Now, let us describe the Nystr\"om method for solving equation \eqref{EqInFrW} which we also write as
$$(I-\mathcal{K})f=g, \qquad  \mathcal{K}=K^1+K^2H$$
where $H$ is given in \eqref{H} and
\[(K^1f)(y)= \int_{-1}^1{ k_1(x,y) f(x)}dx, \qquad   (K^2f)(y)= \int_{-1}^1{ \psi(x)k^*(|x-y|) f(x)}dx. \]
For this purpose, we introduce the  operators
\begin{equation}\label{defK*m}(K_{m}^1f)(y)= \sum_{k=1}^{m}{ \lambda_k k_1(x_k,y) f(x_k)}, \qquad  (K^*_{m}f)(y)= \sum_{k=1}^{m}{ c_{ k}(y) f(x_k)},
\end{equation}
and consider the  equation
\begin{equation*}
(I- \mathcal{K}_m)f_m=g, \qquad \mathcal{K}_m=K^1_m+K^*_mH,
\end{equation*}
i.e.
\begin{equation}\label{fm=G*m}
f_m= \mathcal{G}_mf_m, \quad \mbox{with} \quad \mathcal{G}_mf_m= (\mathcal{K}_mf_m) +g,
\end{equation}
where $f_m$ is an unknown function.

At this point, in order to compute the solution $f_m$ of equation \eqref{fm=G*m} which has the explicit form
\[ f_m(y) - \sum_{k=1}^{m}{\lambda_k k_1(x_k,y) f_m(x_k)}-\sum_{k=1}^{m}{ c_k(y) h(x_k, f_m(x_k))}=g(y),\]
 we collocate it  at the points $x_i,i=1,...,m,$ obtaining the 
 nonlinear system
\begin{equation} \label{NonLinSisDS}
\sum_{k=1}^{m}[\delta_{ik}- \lambda_k k_1(x_k,x_i)]a_k - \sum_{k=1}^{m}{c_k(x_i)  h(x_k, a_k)}= g(x_i), \quad i=1,\ldots,m,
\end{equation}
in the $m$ unknowns $a_i= f_m(x_i)$. \\
The solution $(a_1^*,\ldots,a_m^*)$ allows us to construct the Nystr\"om interpolant as follows
\[ f_m(y)= \sum_{k=1}^{m}{\lambda_k k_1(x_k,y) a^*_k }+ \sum_{k=1}^{m}{ c_k(y)h(x_k,a_k^*)}+g(y).
\]
The difficulty in applying this procedure is the construction of the entries of the matrix of the nonlinear system \eqref{NonLinSisDS} and, in particular, the computation of the constants $c_k(y)$, $k=1,\ldots,m$, for fixed $y$. Indeed, by their definition, the crucial point is the computation of the modified moments  (\ref{modmom}), that can be carried out only for some special kernels $k^*$.
Fortunately, for kernels of type (\ref{k*}) the modified moments can be exactly computed by using the well known recurrence relations for Legendre poloynomials \cite{szego}.
Moreover, as we have already underlined, in this case Theorem \ref{TeoEsUn} still holds true and therefore Theorem \ref{teo:stability} follows also for equation (\ref{fm=G*m}) and this means that it is unisolvent.

Concerning the convergence, in the special case (\ref{k*}) we have the following result.

\begin{theorem}\label{teo:conv*}
	Assume the assumptions of Theorem \ref{teo:stability} let be true.
	Let $f^*$ be the unique fixed point of $\mathcal{G}=K^1+\widetilde{K}^\mu H+g$ in $B(f^*,\delta)$, for some $\delta>0$. Assume that $1$ is not an eigenvalue of $\mathcal{G'}(f^*)$ and in addition for some $\lambda\geq 0$, $-1<\mu \leq 0$,
	\begin{align}
		& \sup_{x \in [-1,1]} \|k_1(x,\cdot)\|_{Z^\lambda}< \infty, \qquad \sup_{y \in [-1,1]} \|k_1(\cdot,y)\|_{Z^\lambda}< \infty,\label{assK1} \\
		& g \in Z^s, \quad \mbox{where} \ s=\left\{ 
			\begin{array}{lc}
				\min{\{\lambda, 1+\mu\}}, & \mu\neq 0,\\ \min{\{\lambda, 1-\epsilon\}}, &\mu=0, \ \epsilon>0 \   \mbox{arbitrarily small.}
			\end{array}
			\right.
			\nonumber
	\end{align}
Then, for $m$ sufficiently large (say $m \geq m_0$), denoted by $f_m$ the unique fixed point of $\mathcal{G}_m$ in $B(f^*,\epsilon)$, with $0 < \epsilon \leq \delta$,  we have
	\begin{equation*} 
		||f^* -f_m||_\infty = \mathcal{O} \left(  \frac{1}{m^s} \right).
	\end{equation*}
\end{theorem}
\begin{remark}
The assumptions of Theorem \ref{teo:conv*} assure that $f^* \in Z^s$, $0<s<1$. If $\mu>0$ and $\lambda>1$, then $s$ could be greater than $1$. In this case, we can apply Theorem \ref{teo:conv} recalling that $f^* \in W^{[s]}$.
\end{remark}
\subsection{Numerical experiments}
In this subsection we give three numerical experiments to show the performance of the method. As the regular case we consider the relative discrete errors as in \eqref{RelDErr}. To solve the numerical system \eqref{NonLinSisDS}, we used the classical Newton method or the Matlab function \verb"fsolve".
\begin{example}\label{Example7}
\rm
Consider the equation proposed \cite{Mandal1}
\[ f(y)-\int_{0}^{1}{|x-y|^{-\frac{1}{2}} (f(x))^2}dx=  g(y),  \quad y \in [0,\,1],
\]
where $g(y)= [y(1-y)]^{\frac{1}{2}} ]+ \frac{16}{15}y^{\frac{5}{2}}+ 2y^2(1-y)^{\frac{1}{2}} + \frac{4}{3}y(1-y)^{\frac{3}{2}}+ \frac{2}{5}(1-y)^{\frac{5}{2}}-\frac{4}{3}y^{\frac{3}{2}}-2 y(1-y)^{\frac{1}{2}}
- \frac{2}{3} (1-y)^{\frac{3}{2}}$ and the exact solution is $f(x)=[x(1-x)]^{\frac{1}{2}}$. The best convergence error in \cite{Mandal1}  is $10^{-4}$ whereas our method produce an error $\mathcal{E}_4=1.97e-14$. The numerical results definitely overcome the theoretical expectation but this is due to the fact that the function $h(x,f(x)=x(1-x)$ is a polynomial and the product rule is exact.
\end{example}

\begin{example}\label{Example8}
\rm
Let us test the method to the following equation already considered in \cite{Mandal1},
\[ f(y)- \int_{0}^{1}{ \log{|x-y|} \sin{( \pi f(x))} }dx=  1, \quad y \in [0,1],
\]
where the exact solution is $f(x)=1$.  
Transformed the equation into $[-1,1]$, it becomes
\begin{align*}
f \left ( \frac{y+1}{2} \right)&- \frac{1}{2} \int_{-1}^{1}{ \log{|y-x|} \sin{\left( \pi f \left(\frac{x+1}{2} \right) \right)}}dx \\ &- \frac{\log{2}}{2} \int_{-1}^{1}{ \sin{ \left( \pi f \left ( \frac{x+1}{2} \right) \right)}}dx=  1, \qquad |y| \leq 1,
\end{align*}
The best convergence  error in \cite{Mandal1} is $10^{-4}$ whereas in our case we have $\mathcal{E}_4=6.66e-16$.
\end{example}

\begin{example}\label{Example9}
\rm
Consider the equation
\[ f(y) - \int_{-1}^{1}{ x^2y f(x)}dx- \int_{-1}^{1}{ |x-y|^{-1/2} \frac{1}{1+ f(x)^2}} dx=  \sqrt{y+1}, \quad \vert y\vert \leq 1.
\]
In this case, we do not known the exact solution. For increasing values of $m$,  Table \ref{table:Example7} reports the relative errors  exhibiting a better performance than the expected one according to the theoretical estimate, which is $\mathcal{O}(m^{-\frac{1}{2}})$, as also confirmed by the estimated order of convergence reported in the last column.

\begin{table}
\caption{The numerical results for Example \ref{Example9} \label{table:Example7}}
\begin{tabular}{llll}
\hline\noalign{\smallskip}
				$m$ &  $\mathcal{E}_m$ & iter & $EOC_m$\\
\noalign{\smallskip}\hline\noalign{\smallskip}
		           8 	& 2.93e-03 	& 6 & 2.01e+00 \\
                  16 	& 7.81e-04 	& 6 & 1.91e+00  \\
                  32 	& 2.03e-04 	& 6 & 1.94e+00  \\
                  64 	& 5.16e-05 	& 6 & 1.98e+00  \\
                  128 	& 1.25e-05 	& 6 & 2.05e+00 	 \\
                  256 	& 2.51e-06 	& 6 & 2.31e+00 	\\
\noalign{\smallskip}\hline
\end{tabular}
\end{table}
\end{example}

\section{An application to Boundary Integral Equations}\label{sec:BIE}
In this section, we show an application of the Nystr\"om method described in Section \ref{sec:weakly} for the numerical solution of a nonlinear boundary integral equation (BIE) arising from the reformulation of a nonlinear boundary value problem (BVP) for Laplace's equation.

Let us consider the interior Neumann problem over a bounded simply connected planar domain $D \subset \RR^2$ with smooth boundary $\Gamma$. It consists in finding a function
$u \in C^2(D) \cap C^1(\overline{D})$ that satisfies
\begin{equation} \label{BVP}
\begin{cases}
\Delta u(P)=0, & \qquad P \in D, \vspace{0.1cm} \\
\displaystyle \frac{\partial u(P)}{\partial n_P}=-\bar{h}(P,u(P))+\bar{g}(P), & \qquad P \in \Gamma,
\end{cases}
\end{equation}
where  $n_P$ denotes the exterior unit normal to $\Gamma$ at the point $P$, while the function $\bar{h}(P,v)$ defined in $\Gamma \times \RR$ is  nonlinear in $v$ and is assumed sufficiently smooth.

It is known that (see, for instance, \cite{AtkChan}) the harmonic function $u$ satisfying \eqref{BVP} is the solution of the following nonlinear BIE of the second kind
\begin{equation}\label{BIE}
\begin{split}
u(P)-\frac{1}{\pi}\int_{\Gamma}u(Q)\frac{\partial}{\partial n_Q}[\log{|P-Q|}]&d\sigma(Q) -\frac{1}{\pi}\int_{\Gamma}\bar{h}(Q,u(Q))\log{|P-Q|}d\sigma(Q)\\
&=\frac{1}{\pi}\int_{\Gamma}\bar{g}(Q)\log{|P-Q|}d\sigma(Q), \qquad P \in \Gamma
\end{split}
\end{equation}
which can be deduced from  Green's representation formula for $u$
\begin{equation}\label{GreenFormula}
\begin{split}
u(P)=\frac{1}{2 \pi}\int_{\Gamma}u(Q)\frac{\partial}{\partial n_Q}&[\log{|P-Q|}]d\sigma(Q) \\ &-\frac{1}{2 \pi}\int_{\Gamma}\frac{\partial u(Q)}{\partial n_Q}\log{|P-Q|}d\sigma(Q), \quad P \in D,
\end{split}
\end{equation}
taking into account the boundary condition in \eqref{BVP}. \newline
Once equation \eqref{BIE} have been solved, one can use the known function $u$ on $\Gamma$ along with its known normal derivative on the boundary given in \eqref{BVP} in order to compute the unknown solution $u$ on the domain $D$ by means of formula \eqref{GreenFormula}.

Hence, we are interested in the numerical solution of \eqref{BIE}. In order to transform the BIE \eqref{BIE} into an equivalent 1D integral equation on the interval $[-1,1]$,  firstly we introduce a parametric representation of the curve $\Gamma$  
\begin{equation} \label{parameterization}
\boldsymbol{\gamma}(x)=(\xi(x),\eta(x)) \in \Gamma, \qquad x \in [-1,1].
\end{equation}
We assume that $\gamma$ traverses $\Gamma$ in a counter-clockwise direction (i.e. it is such that the domain $D$ is on the left of $\Gamma$) and $\xi,\eta \in C^2[-1,1]$, with
\begin{equation*}
|\gamma'(x)| \neq 0, \qquad \forall x \in [-1,1].
\end{equation*}
Moreover, in order to achieve higher orders of convergence for the numerical method and, consequently, more accurate approximations of the solution, we adopt some already known regularization strategies (see, for instance, \cite{ElJeSlSt,MS,Scuderi,FermoLaurita2,Laurita}) considering a smoothing
transformation $\phi(x)$, such that
\begin{equation} \label{smoothing}
\phi(x)=\left\{
\begin{array}{ll}
-1+(x+1)^q, & \quad x \in [-1,-1+\epsilon],\\
1-(1-x)^q, & \quad  x \in [1-\epsilon,1],\\
\end{array} \right.
\end{equation}
for some small $\epsilon>0$ and some smoothing exponent $q \geq 1$.
Note that in the case $q=1$ we have $\phi(x)=x$, which means that no smoothing transformation is applied.

Then, by introducing in \eqref{BIE} the change of variables $x={\bar \gamma}(x)$ and $y={\bar \gamma}(y)$ with
\begin{equation} \label{smoothparametrization}
{\bar \gamma}(x)=\gamma(\phi(x))=(\xi(\phi(x)),\eta(\phi(x)))=:(\bar{\xi}(x),\bar{\eta}(x)) \qquad x \in [-1,1].
\end{equation}
we can rewrite the BIE as follows
\begin{equation} \label{1DIntEq}
\begin{split}
u(\bar{\gamma}(y))&-\int_{-1}^1 k_1(x,y) u(\bar{\gamma}(x))dx \\
&-\frac{1}{\pi}\int_{-1}^1\bar{h}(\bar{\gamma}(x),u(\bar{\gamma}(x)))|\bar{\gamma}'(x)|\log{|\bar{\gamma}(y)-\bar{\gamma}(x)|}dx \\
&=-\frac{1}{\pi}\int_{-1}^1\bar{g}(\bar{\gamma}(x))|\bar{\gamma}'(x)|\log{|\bar{\gamma}(y)-\bar{\gamma}(x)|}dx,
\end{split}
\end{equation}
where
\begin{equation*} \label{kernelk1}
k_1(x,y)=\frac{1}{\pi}\left\{
\begin{array}{ll}
\ds\frac{\bar{\eta}'(x)[\bar{\xi}(x)-\bar{\xi}(y)]-\bar{\xi}'(x)[\bar{\eta}(x)-\bar{\eta}(y)]}
        {[\bar{\xi}(x)-\bar{\xi}(y)]^2+[\bar{\eta}(x)-\bar{\eta}(y)]^2}, & \qquad x \neq y, \vspace*{0.5em}\\
\ds\frac{1}{\pi}\frac{\bar{\xi}'(x)\bar{\eta}''(x)-\bar{\eta}'(x)\bar{\xi}''(x)}{2[\bar{\xi}'(x)^2+\bar{\eta}'(x)^2]}, &\qquad x=y.
\end{array} \right.
\end{equation*}
Now, first we split the logarithmic kernel $\log{|\bar{\gamma}(y)-\bar{\gamma}(x)|}$ as follows
\begin{equation}\label{log}
\log{|\bar{\gamma}(y)-\bar{\gamma}(x)|}=\log{\frac{|\bar{\gamma}(y)-\bar{\gamma}(x)|}{|x-y|}}+\log{|x-y|}.
\end{equation}
Then, following a numerical trick in \cite{MonegatoScuderi} (see, also, \cite{FermoLaurita}) in order to avoid numerical cancellation,  when $|x-y|<eps$ ($eps$ denotes the machine precision) we use the approximation
\begin{equation*}
\log{\frac{|\bar{\gamma}(y)-\bar{\gamma}(x)|}{|x-y|}}\simeq \log{|\bar{\gamma}'(x)|}.
\end{equation*}
Now, setting $f(x)=u(\bar{\gamma}(x))$, $h(x,f(x))=\bar{h}(\bar{\gamma}(x),u(\bar{\gamma}((x))$,
\begin{equation*} \label{kernell}
\rho(x,y)=\left\{
\begin{array}{ll}
\ds \frac{1}{\pi}|\bar{\gamma}'(x)|\log{|\bar{\gamma}'(x)|} & \qquad |x-y|<eps,\vspace*{0.5em}\\
\ds \frac{1}{\pi} |\bar{\gamma}'(x)|\log{\frac{|\bar{\gamma}(y)-\bar{\gamma}(x)|}{|x-y|}}, & \qquad \mathrm{otherwise},
\end{array} \right.
\end{equation*}
$\psi(x)=\frac{1}{\pi}|\bar{\gamma}'(x)|$, $k_2(x,y)=\rho(x,y)+\psi(x)\log{|x-y|}$, and, finally,
\begin{equation} \label{rhsideBIE}
g(y)=-\frac{1}{\pi}\int_{-1}^1\bar{g}(\bar{\gamma}(x))|\bar{\gamma}'(x)|\log{|\bar{\gamma}(y)-\bar{\gamma}(x)|}dx,
\end{equation}
the integral equation \eqref{1DIntEq} takes the form  \begin{equation} \label{IntEqBIE}
f(y)=(\mathcal{G}f)(y), \qquad \mathrm{with}  \quad (\mathcal{G}f)(y)=(K_1f)(y)+(K_2Hf)(y)+g(y)
\end{equation}
with the operators $K_i$, $i=1,2$, and $H$ defined as in \eqref{K} and \eqref{H}, respectively.

In order to approximate the solution of \eqref{IntEqBIE}, we apply a Nystr\"om type method which is a combination of the methods described in sections \ref{sec:smooth} and \ref{sec:weakly}. More precisely, it consists in solving the following approximating equation
\begin{equation} \label{DiscreteIntEqBIE}
f_m(y)=(\mathcal{G}_mf_m)(y)
\end{equation}
with the operator $\mathcal{G}_m$ defined as
\begin{equation*}
(\mathcal{G}_mf)(y)=(K_m^1f)(y)+(K_m^2Hf)(y)+g(y),
\end{equation*}
where $K_m^1$ is the operator given in \eqref{Km}, while $K_m^2$ is the operator defined as follows
\begin{equation*}
(K_m^2f)(y)=\sum_{k=1}^m\left[\lambda_k \rho(x_k,y)+c_k(y)\psi(x_k)\right]f(x_k).
\end{equation*}
The collocation of equation \eqref{DiscreteIntEqBIE} at the Legendre zeros leads to the nonlinear system
\begin{equation} \label{nonlinearsystemBIE}
f_m(x_k)=(\mathcal{G}_mf_m)(x_k), \qquad k=1,\ldots,m,
\end{equation}
whose solutions are
$f_m(x_i)$, $i=1,\ldots,m$. \newline
Once computed, these values can be used in order to construct an approximation of the harmonic function $u$, solution of the boundary value problem \eqref{BVP}, at any point $P$ of the interior domain $D$.  First, using the parameterization \eqref{smoothparametrization} in \eqref{GreenFormula}, for any $P\equiv(x_P,y_P)\in D$ we represent the potential $u(P)$ as
\begin{equation} \label{DoublePotential}
\begin{split}
u(x_P,y_P)&=\frac{1}{2}\int_{-1}^1\frac{\bar{\eta}'(x)[\bar{\xi}(x)-x_P]-\bar{\xi}'(x)[\bar{\eta}(x)-y_P]}
        {[\bar{\xi}(x)-x_P]^2+[\bar{\eta}(x)-y_P]^2}f(x)dx\\
        &+ \frac{1}{2}\int_{-1}^1 [h(x,f(x))-\bar{g}(\bar{\gamma}(x))]|\bar{\gamma}'(x)|\log{\left|(x_P,y_P)-\bar{\gamma}(x)\right|}dx,
\end{split}
\end{equation}
where $f$ is the solution of the integral equation \eqref{IntEqBIE}.
Then, we approximate the potential $u$ in \eqref{DoublePotential}, by the following function
\begin{equation} \label{DiscreteDoublePotential}
\begin{split}
&u_m(x_P,y_P) \\ &=\frac{1}{2}\sum_{k=1}^m \lambda_k\frac{\bar{\eta}'(x_k)[\bar{\xi}(x_k)-x_P]-\bar{\xi}'(x_k)[\bar{\eta}(x_k)-y_P]}
        {[\bar{\xi}(x_k)-x_P]^2+[\bar{\eta}(x_k)-y_P]^2}f_m(x_k)\\
        &+ \frac{1}{2}\sum_{k=1}^m \lambda_k  \left[h(x,f_m(x_k))-\bar{g}(\bar{\gamma}(x_k))\right]|\bar{\gamma}'(x_k)|\log{\left|(x_P,y_P)-\left(\bar{\xi}(x_k),\bar{\eta}(x_k)\right)\right|},
\end{split}
\end{equation}
where the solutions $f_m(x_i)$, $i=1,\ldots,m$, of the solved nonlinear system are employed. \newline
We observe that in the computation of the right-hand sides of such system the function $g(y)$ given in \eqref{rhsideBIE} needs to be evaluated at the collocation points. When we are not able to compute analytically the integral, proceeding as in \cite{FermoLaurita}, we approximate it taking into account \eqref{log} and  using a proper combination of the Gauss-Legendre formula and a product quadrature rule with a large number of knots.

\subsection{Numerical experiments}
In this subsection, we are going to show some numerical examples in which the method described in sections \ref{sec:smooth}, \ref{sec:weakly} has been applied for approximating the solution of the interior Neumann problem \eqref{BVP} in some planar domain $D$ with smooth boundary $\Gamma$.

The reported error $\|u-u_m\|_{\Gamma}$ is the maximum error at the node points on $\Gamma$ while
the error $\|u-u_m\|_D$ represents the maximum error at 600 points sampled randomly in the interior domain $D$.
Moreover, in our tests, whenever necessary, we have used as smoothing transformation $\phi(x)$ (see \eqref{smoothing}) the following one adopted in \cite{MS,Scuderi}
\begin{equation*} \label{phiq}
\phi(x)=\frac{2\int_{-1}^x \left(1-t^2\right)^{q-1}dt}{\int_{-1}^1 \left(1-t^2\right)^{q-1}dt}-1, \qquad   x \in [-1,1], \quad q \geq 1.
\end{equation*}

\begin{example} \label{ExampleBIE1}
\rm
We consider the problem \eqref{BVP} defined on the planar region $D$ bounded by the ellipse $\Gamma$ of equation
\begin{equation*}
\frac{x^2}{a^2}+\frac{y^2}{b^2}=1
\end{equation*}
for given values of $(a,b)$. We choose the function $\bar{h}$ as follows  (see \cite{AtkChan})
\begin{equation*}
\bar{h}(P,v)=v+\sin{v},
\end{equation*}
and the function $\bar{g}$ such that the exact solution of \eqref{BVP} is the harmonic function
\begin{equation} \label{exactSolBIE}
u(x,y)=e^x\cos{y}.
\end{equation}
Table \ref{table1BIE} contains the numerical results obtained by applying the numerical method \eqref{nonlinearsystemBIE}, \eqref{DiscreteDoublePotential} in the case $(a,b)=(1,2)$.

\begin{table}
\caption{Errors for the potential $u$ in Example \ref{ExampleBIE1}}
\label{table1BIE}
\begin{tabular}{lllll}
\hline\noalign{\smallskip}
	 $ $     &\multicolumn{2}{c}{$q=1$}    &\multicolumn{2}{c}{$q=2$}  \\ \hline 
     $m$ &  $\|u-u_m\|_{\Gamma}$ & $\|u-u_m\|_D$ & $\|u-u_m\|_{\Gamma}$ & $\|u-u_m\|_D$ \\ 
\noalign{\smallskip}\hline\noalign{\smallskip}
8   &	 6.93e-02  &  2.71e-01  & 4.59e-01   & 3.58e-01  \\
            16  & 	 2.36e-03  &  5.94e-02  & 1.14e-02   & 1.42e-01 \\
            32  &	 3.94e-04  &  4.98e-03  & 6.60e-05   & 1.48e-02 \\
            64  &	 1.01e-04  &  3.65e-05  & 5.92e-07   & 1.21e-03  \\
            128 & 	 2.56e-05  &  3.89e-08  & 3.76e-08	 & 2.59e-05 \\
            256 &	 6.44e-06  &  2.45e-09  & 3.19e-09	 & 2.84e-09 \\
            512 &	 1.61e-06  &  1.53e-10  & 1.98e-09	 & 1.71e-14 \\
					           
\noalign{\smallskip}\hline
\end{tabular}
\end{table}
\end{example}

\begin{example}      \label{ExampleBIE2}
\rm
In this second example, again taken from  \cite{AtkChan}, we solve the interior Neumann problem \eqref{BVP} defined on the same domain $D$ considered in the previous example. Here, we assume but assuming the boundary condition defined by the nonlinear function 
\begin{equation*}
\bar{h}(P,v)=|v|v^3
\end{equation*}
and the function $\bar{g}$ chosen such that the exact solution $u$ of \eqref{BVP} is the one given in \eqref{exactSolBIE}.\newline
The obtained numerical results are shown in Table \ref{table2BIE}.

\begin{table}
\caption{Errors for the potential $u$ in Example \ref{ExampleBIE2}}\label{table2BIE}
\begin{tabular}{lllll}
\hline\noalign{\smallskip}
	 $ $     &\multicolumn{2}{c}{$q=1$}    &\multicolumn{2}{c}{$q=2$}  \\ \hline 
     $m$ &  $\|u-u_m\|_{\Gamma}$ & $\|u-u_m\|_D$ & $\|u-u_m\|_{\Gamma}$ & $\|u-u_m\|_D$ \\ 
\noalign{\smallskip}\hline\noalign{\smallskip}
         8 	&  5.27e-01  &	 6.28e-01 &  6.71e-01   & 8.71e-01       \\
         16 &  5.65e-02  &	 9.57e-02 &  9.83e-02  & 1.95e-01      \\
         32 &  1.23e-03  &	 6.01e-03 &  1.84e-03  & 1.52e-02      \\
         64 &  7.58e-04  &	 2.23e-04 &  3.38e-05  & 1.20e-03       \\
         128&  3.16e-04  &	 1.00e-06 &  6.53e-07  & 2.60e-05      \\
         256&  9.89e-05  &	 4.54e-08 &  4.21e-08  & 5.54e-09      \\
         512&  2.70e-05  &	 3.06e-09 &  1.53e-09  & 3.18e-10     \\
\noalign{\smallskip}\hline
\end{tabular}
\end{table}
\end{example}

\begin{example}   \label{ExampleBIE3}
\rm
We consider the  amoeba-like domain $D$ bounded by the curve $\Gamma$ having the following parametric representation
\begin{equation*}
\gamma(x)=R(\pi(x+1))e^{{\mathrm i}\pi(x+1)}, \qquad x \in [-1,1]
\end{equation*}
with $R(x)=e^{\cos{x}}\cos^2{2x}+e^{\sin{x}}\sin^2{2x}$ (see Figure \ref{Fig:amoeba}).

\begin{figure*}
\includegraphics[width=0.8\textwidth]{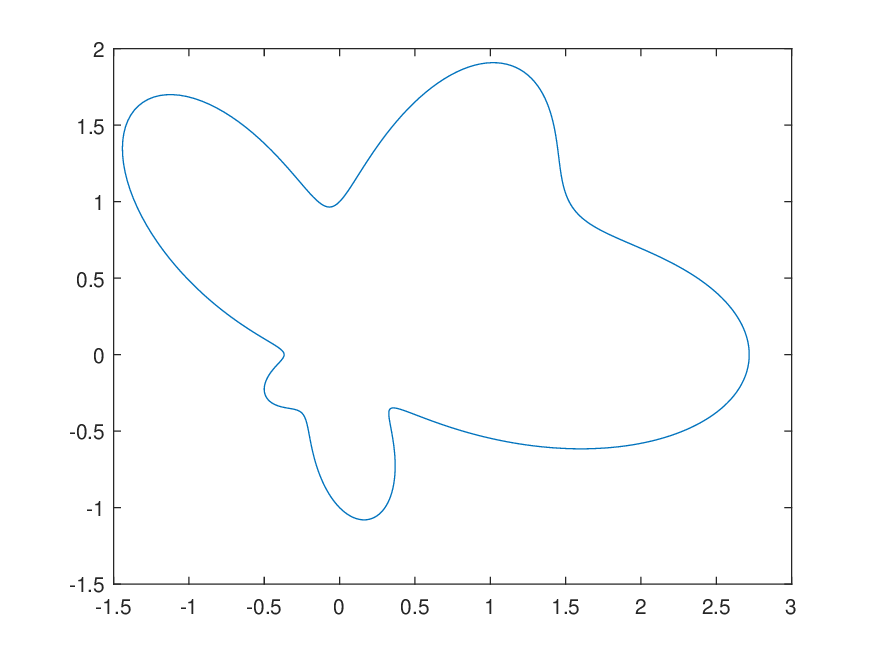}
\caption{The boundary $\Gamma$ of the domain $D$ in Example \ref{ExampleBIE3} \label{Fig:amoeba}}
\end{figure*}

We assume as exact solution of the BVP \eqref{BVP}, with the nonlinear function $\bar{h}(P,v)=v^3$ ,
the function
\begin{equation*}
u(x,y)=\sin(x)\cosh(y)
\end{equation*}
and determine the corresponding function $\bar{g}$ at the right-hand side of the Neumann boundary condition. The errors occurred in the computation of the solution at the nodes of curve $\Gamma$ and at some random points in the interior domains are reported in Table \ref{table3BIE}.

\begin{table}
\caption{Errors for the potential $u$ in Example \ref{ExampleBIE3}}
\label{table3BIE}
\begin{tabular}{llll}
\hline\noalign{\smallskip}
 $ $     &\multicolumn{2}{c}{$q=1$} \\ 	\hline
				$m$ &  $\|u-u_m\|_{\Gamma}$ & $\|u-u_m\|_D$ \\                      
\noalign{\smallskip}\hline\noalign{\smallskip}
 16 	 &	     4.23e-01  &	 4.72e-01 \\
			 32 	 & 	 1.99e-01  &	 1.09e-01 \\
			 64 	 & 	 5.68e-03  &	 8.61e-03 \\
			128   &  	 3.25e-05  &	 1.38e-04 \\
			256 	 & 	 1.94e-07  &	 6.71e-08 \\
			512 	 & 	 4.86e-08  &	 1.42e-12 \\
\noalign{\smallskip}\hline
\end{tabular}
\end{table}
 
We observe that we have obtained accurate numerical results already without using any smoothing transformation ($q=1$), even better that those obtained for $q>1$, when one needs larger values of $m$ in order to achieve higher accuracy.
\end{example}

\section{Proofs}\label{sec:proofs}
In this section, we collect the proofs of all our main results.
From now on, for the sake of simplicity, when we handle with the norm of an operator $T$,  we will omit the subscript $C^0 \to C^0$
in the norm, that is we set $\|T\|=\|T\|_{C^0 \to C^0}$.

\begin{proof}\emph{of Theorem \ref{teo:K}.} 
For any fixed $x \in [-1,1]$, let $P_m(x, y)$ and $Q_m(x, y)$ be the polynomial of best approximation with respect to the variable $y$ of $k_1(x, y)$ and $k_2(x, y)$, respectively, i.e.
$$E_m(k_1(x,\cdot))=\|k_1(x,\cdot) -P_m(x,\cdot)\|_\infty, \quad E_m(k_2(x,\cdot))=\|k_2(x,\cdot) -Q_m(x,\cdot)\|_\infty,$$
and introduce the operator $\mathcal{\tilde{K}}_{m}$ defined as follows
\begin{equation*}
(\mathcal{\tilde{K}}_{m}f)(y)= \int_{-1}^{1}{P_m(x,y) f(x)dx}+\int_{-1}^{1}{Q_m(x,y) h(x,f(x))} dx.
\end{equation*}
Denoted by $\mathcal{K}=K^1+K^2H$ where $K^i$, $i=1,2$, are given in \eqref{K} and $H$ is defined in \eqref{H}, we have
\begin{align*}
&|(\mathcal{K}f)(y)- (\mathcal{\tilde{K}}_{m} f)(y)| \\
&\leq \int_{-1}^{1}{|k_1(x,y)- P_m(x,y)|}dx +
\int_{-1}^{1}{|k_2(x,y)- Q_m(x,y)| \,  |h(x,f(x))|}dx   \\
& \leq 2  \sup_{|x| \leq 1}{|k_1(x,y)-P_m(x,y)|}+ \sup_{|x| \leq 1}{{|k_2(x,y)-Q_m(x,y)|} \int_{-1}^{1}{|h(x,f(x))|}dx}.
\end{align*}

Then, since $\mathcal{\tilde{K}}_{m}f$ is a polynomial of degree at most $m$, for any $f \in C^0$, it follows
\begin{equation*}
E_m(\mathcal{K}f) \leq 2\sup_{|x| \leq 1}{E_m(k_1(x,\cdot))} +\sup_{|x| \leq 1}{E_m(k_2(x,\cdot))}  \int_{-1}^{1}{|h(x,f(x))|}dx.
\end{equation*}
Therefore, by the assumptions on the kernels $k_i$, $i=1,2$, and on the function $h$, and by applying \eqref{Ch.Sobolev}, we deduce that $ \mathcal{K}f \in W^r$. Hence, if $g \in W^r$, the solution of equation \eqref{Hamm_op} $f=\mathcal{K}f+g \in W^r$.
\end{proof}

\begin{proof}\emph{of Theorem {\ref{theoWeak}}.}
	For simplicity set
	\[
	\mathcal{G}f=K^1f+K^2Hf+g,
	\]
	with $K^1$ and $K^2$ defined as in (\ref{K}), and with $k_2(x,y)=\psi(x)k^*(|x-y|)$, while $H$ is defined in (\ref{H}).  Then the solvability of equation (\ref{EqInFrW}) is equivalent to the existence of a fixed point of the operator $\mathcal{G}$ in $C^0$.
	
Then, for $f_1 $, $f_2 \in C^0$.  We get
	\begin{eqnarray*}
		\|	\mathcal{G}f_1-	\mathcal{G}f_2\|_\infty &\leq&  \max_{y\in[-1,1]} \left[ \int_{-1}^1 |k_1(x,y)||f_1(x)-f_2(x)|\ dx \right.\\ 
		&+ &\left. \int_{-1}^1 |\psi(x) k^*(|x-y|)||h(x,f_1(x))-h(x,f_2(x))|\ dx\right] \\
		& \leq & \|f_1-f_2\|_\infty \max_{y\in[-1,1]}\left( \int_{-1}^1 |k_1(x,y)|\ dx  \right. \\  &+& \left. \sup_{x\in[-1,1]}\|h_v(x,\cdot)\|_\infty \int_{-1}^1 \psi(x)|k^*(|x-y|)|\ dx\right)
		\end{eqnarray*}
		Therefore, under the assumption \textbf{[H4]}, we deduce that $	\mathcal{G}$ is a contraction mapping on $C^0$ and consequently it has a unique fixed point.
	\end{proof}
	
\begin{proof}\emph{of Theorem {\ref{classeDebSing}}.}
We want to estimate $E_m(\mathcal{K}f)$ in order to apply (\ref{Ch.Zyg}) and deduce the class of the function $\mathcal{K}f$.

First of all we underline that 
\[
E_m(\mathcal{K}f) \leq E_m(K^1f)+ E_m(\widetilde{K}^\mu Hf).
\]
About $E_m(K^1f)$, proceeding, for instance, as done in the proof of Theorem \ref{teo:K} we get
\[
E_m(K^1f) \leq \mathcal{C} \sup_{|x|\leq 1} E_m(k_1(x,\cdot)), \quad \mathcal{C}\neq\mathcal{C}(m,f)
\]
and hence, under the assumptions \eqref{assk1}, we get
\[E_m(K^1f)\leq  \frac{\mathcal{C}}{m^\lambda}, \quad \mathcal{C}\neq \mathcal{C}(m,f).\]
Therefore from (\ref{Ch.Zyg}) we deduce $K^1f \in Z^\lambda$.

Consider now $ E_m(\widetilde{K}^\mu Hf)$. Using inequality (\ref{jackson}) we get
\begin{equation}\label{stimaOmega}
 E_m(\widetilde{K}^\mu Hf) \leq \mathcal{C} \int_0^{\frac 1 m} \frac{\Omega^k_\varphi(\widetilde{K}^\mu Hf,t)}{t} \ dt, \quad \mathcal{C}\neq \mathcal{C}(m,f), \quad k \geq 1.
\end{equation}
Thus we have to estimate the main part of the modulus of continuity of $\widetilde{K}^\mu Hf$. We can proceed following step by step the proof of Lemma 4.1 in \cite{MRT}.
Therefore, since we are assuming that \textbf{[H1]} holds true and $\psi \in C^0$, we get
\[
\Omega_\varphi^k(\widetilde{K}^\mu Hf,t) \leq \mathcal{C} \|h(\cdot,f)\|_\infty \|\psi\|_\infty\left\{\begin{array}{lc} t^{1+\mu}, & \mu \neq 0, \  k>1+\mu,\\ 
	t \log t^{-1}, & \mu=0, \  k\geq 1,
\end{array}\right.\]
and, consequently, by (\ref{stimaOmega}) we have
\[
 E_m(\widetilde{K}^\mu Hf) \leq \mathcal{C}\left\{\begin{array}{lc} \displaystyle \frac{1}{m^{1+\mu}}, & \ \mu \neq 0, \  k>1+\mu,\\  & \\
 	\frac{\log m}{m}, & \mu=0, \ k\geq 1.
 \end{array}\right.
\]
From these estimates and (\ref{Ch.Zyg}) we get that $\widetilde{K}^\mu Hf \in Z^r$, where $r=1+\mu$ when $\mu \neq 0$, while $r=1-\epsilon$, with $\epsilon$ sufficiently small, when $\mu=0$, and the theorem follows.
\end{proof}
	
	\begin{proof}\emph{of Theorem {\ref{teo:K_Km}}.} 
		First, note that the well-known Gauss-Legendre formula \eqref{gauss} is convergent. Therefore, the sequences $\{K^i_m\}_m$ $i=1,2$ given in \eqref{Km} are collectively compact and pointwise convergent to the integral operators $K^i$  $i=1,2$ defined in \eqref{K}; see for instance \cite[Theorem 12.8]{Kress}. Consequently, the assertion follows taking into account that the operator $\mathcal{K}_m$ is the sum of the operator $K^1_m$ and the composition of $K^2_m$ with the operator $H$ in \eqref{H}  (see, for instance, \cite[p.74]{Krasnoseski}).
	\end{proof}
\begin{proof}\emph{of Theorem {\ref{teo:stability}}.} 
The existence of a fixed point $f_m \in B(f^*,\epsilon)$ for the operator $\mathcal{G}_m$,  for $m$ sufficiently large, say $m \geq m_0$, follows by \cite[Theorem 3]{KEAtkinson}. In fact, the first four hypothesis of Theorem 3 in \cite{KEAtkinson} are guaranteed  by Theorem \ref{teo:K_Km} and the use of our assumptions.
In addition, the fixed point $f^*$ has nonzero index and is isolated, i.e. $f^* \in B(f^*,\epsilon_0)$ with  $0<\epsilon_0 \leq \delta$. This can be deduced by \cite[Theorem 21.6, p.108]{Krasnoseski} or \cite[p.136]{Krasnoseski64}, taking into account that $g$ verifies ${\bf{[G1]}}$, $h$ satisfies ${\bf{[H1]}}$ and ${\bf{[H2]}}$,  and $1$ is not an eigenvalue of $\mathcal{G}'f^*$.

Let us now prove the uniqueness by showing that
\begin{equation}\label{unicity}
\|f-\mathcal{G}_m(f)\|_{\infty}>0, \qquad  \forall f \in B(f^*,\epsilon),
\end{equation}
i.e. $f$ cannot be a solution of \eqref{Gmfm}.

First, let us note that for sufficiently large $m$, the operators $(I-\mathcal{G}'_mf)^{-1}$, $f \in B(f^*,\epsilon)$, $0 <\epsilon \leq \epsilon_0$, exist and are uniformly bounded w.r.t. $m$ , i.e. there exists a constant $\C\neq \C(m)$ such that
\begin{equation}\label{ii}
\|(I-\mathcal{G}'_mf)^{-1}\| < 2\mathcal{C}, \qquad f \in B(f^*,\epsilon), \qquad  m \geq m_0,
\end{equation}
where, for $\phi \in C^0$,
\begin{align*}
[(\mathcal{G}_m'f)\phi](y)& =\sum_{k=1}^m \lambda_k k_1(x_k,y) \phi(x_k) \\ & +\sum_{k=1}^m \lambda_k {k_2(x_k,y)h_{v}(x_k,f(x_k))\phi(x_k)}, \quad y\in [-1,1].
\end{align*}
This can be obtained by proceeding as in the proof of \cite[Theorem 4]{KEAtkinson}   by virtue of [{\bf{H1}}] [{\bf{H2}}] and [{\bf{H3}}].

Now, for $f \in B(f^*,\epsilon)$ we have
$$f-\mathcal{G}_m(f) = [I-\mathcal{G}'_mf_m](f-f_m)-[\mathcal{G}_m(f)-\mathcal{G}_m(f_m)-\mathcal{G}'_mf_m(f-f_m)],$$
and then, using the reverse triangle inequality,
\begin{equation}\label{f-Gm}
\|f-\mathcal{G}_m(f)\|_\infty \geq \|I-\mathcal{G}'_mf_m\| \|f-f_m\|_\infty-\|\mathcal{G}_m(f)-\mathcal{G}_m(f_m)-\mathcal{G}'_mf_m(f-f_m)\|_\infty.
\end{equation}
Therefore, by \eqref{ii}, since $f_m \in B(f^*,\epsilon)$ we have
\begin{align}\label{ii2}
\|I-\mathcal{G}'_mf_m\|  \geq \frac{1}{2 \C}.
\end{align}
In addition, from  {[\bf{H3}]} we can deduce that also $\mathcal{G}''_m$ admits the second Frechet derivative given by
\begin{equation*}
[(\mathcal{G}_m''f)(\phi_1,\phi_2)](y)=  \sum_{k=1}^m \lambda_k k_2(x_k,y)h_{vv}(x_k,f(x_k))\phi_1(x_k)\phi_2(x_k),   \, \phi_1,\phi_2 \in C^0,
\end{equation*}
and we have
$$\max\{\|\mathcal{G}''f\|,\|\mathcal{G}_m''f\| \} \leq \mathcal{C}_1, \qquad \mathcal{C}_1=\mathcal{C}_1(f^*,\delta), \quad  \mathcal{C}_1\neq\mathcal{C}_1(m), \quad f \in B(f^*,\delta),$$
with $\| \cdot\|$ denoting the norm of a bilinear form from $C^0 \times C^0 \to C^0$, and $\mathcal{G}''f$  given in \eqref{G''}.

Standard arguments \cite[Chapter 17, p. 500]{KantarovicAkilov} lead to
\begin{align*}
-\|\mathcal{G}_m(f)-\mathcal{G}_m(f_m)-\mathcal{G}'_m f_m(f-f_m)\|_\infty  &\geq -\frac{\mathcal{C}_1}{2} \|f-f_m\|^2_\infty  \\ & \geq - \epsilon \, \mathcal{C}_1 \|f-f_m\|_\infty.
\end{align*}
Therefore, by applying \eqref{ii2} and the above inequality in \eqref{f-Gm}, by fixing $\epsilon <(2 \mathcal{C} \mathcal{C}_1)^{-1}$, we have
\begin{align*}
\|f-\mathcal{G}_m(f)\|_\infty \geq \left[\frac{1}{2 \mathcal{C}}-\epsilon \mathcal{C}_1 \right] \|f-f_m\|_\infty, \qquad m > m_0\qquad  f \in B(f^*,\epsilon),
\end{align*}
namely \eqref{unicity}.
\end{proof}

\begin{proof}\emph{of Theorem \ref{teo:conv}.}  
First let us prove that
\begin{equation}\label{f-fm}
\|f^* -f_m\|_\infty \leq  \mathcal{C} \,  \|(\mathcal{G}-\mathcal{G}_m)f^*\|_\infty,
\end{equation}
where $\C \neq \C(m)$. To this end, let us proceed as in the proof of Theorem 4 in \cite{KEAtkinson}. Then, by \eqref{EqInFrOp} and \eqref{Gmfm} we write
$$f^*-f_m=\mathcal{G}f^*-\mathcal{G}_mf_m.$$
Therefore,
$$(I-\mathcal{G}'_mf^*)(f^*-f_m)=(\mathcal{G}-\mathcal{G}_m)f^*-[\mathcal{G}_m(f_m-f^*)-\mathcal{G}'_m f^*(f_m-f^*)],$$
or equivalently
$$f^*-f_m=(I-\mathcal{G}'_mf^*)^{-1}\{(\mathcal{G}-\mathcal{G}_m)f^*-[\mathcal{G}_m(f_m-f^*)-\mathcal{G}'_m f^*(f_m-f^*)]\},$$
from which it follows
\begin{align*}
\|f^*-f_m\|_\infty & \leq \|(I-\mathcal{G}'_mf^*)^{-1}\| \{\|(\mathcal{G}-\mathcal{G}_m)f^*\|_\infty  \\ &+ \|\mathcal{G}_m(f_m-f^*)-\mathcal{G}'_mf^*(f_m-f^*)\|_\infty\}.
\end{align*}
By applying
$$\|(I-\mathcal{G}'_mf^*)^{-1}\| < \mathcal{C},  \qquad  m \geq m_0$$
and  \cite[Chapter 17, p. 500]{KantarovicAkilov}
$$\|\mathcal{G}_m(f_m-f^*)-\mathcal{G}'_mf^*(f_m-f^*)\|_\infty \leq \frac{\mathcal{C}_1}{2} \|f^*-f_m\|^2_\infty, $$
we have
\begin{align*}
\|f^*-f_m\|_\infty  & \leq  \frac{ \C}{1-\frac{\C \mathcal{C}_1}{2} \|f^*-f_m\|_\infty} \|(\mathcal{G}-\mathcal{G}_m)f^*\|_\infty.
\end{align*}

Taking $\epsilon< (2\C \C_1)^{-1}$, since $f_m \in B(f^*,\epsilon)$, we get
\begin{align*}
\|f^*-f_m\|_\infty  & \leq   \frac{ \C}{1-\frac{\C \mathcal{C}_1}{2} \epsilon} \|(\mathcal{G}-\mathcal{G}_m)f^*\|_\infty
\end{align*}
from which we deduce \eqref{f-fm} being the denominator less than $1/2$.

Now, by using the definition of $\mathcal{G}$ and $\mathcal{G}_m$ and by applying \eqref{em}, we deduce
\begin{align*}
\|f^* -f_m\|_\infty &\leq  \mathcal{C} \left[ \|({K}^1-{K}^1_m)f^*\|_\infty+  \|({K}^2- {K}_m^2)f^*\|_\infty \right] \\ & = \mathcal{C} \left[ \sup_{|y| \leq 1} |e_m(k_1(\cdot,y) f^*)|+\sup_{|y| \leq 1} |e_m(k_2(\cdot,y) h(\cdot,f^*))| \right] \\ & \leq \C \left[\sup_{|y| \leq 1} E_{2m-1}(k_1(\cdot,y) f^*)+
 \sup_{|y| \leq 1} E_{2m-1}(k_2(\cdot,y) h(\cdot,f^*)) \right].
\end{align*}
Therefore, by exploiting the following estimate
$$E_{2m}(f_1 f_2) \leq \|f_1\|_\infty E_m(f_2)+2 \|f_2\|_\infty E_m(f_1), \quad \forall f_1,f_2 \in C^0,$$
we have
\begin{align*}
\|f^* -&f_m\|_\infty \leq  \mathcal{C} \left[ \| f^*\|_\infty \, \sup_{|y| \leq 1} E_{m-1}(k_1(\cdot,y))  + 2 \sup_{|y| \leq 1} \|k_1(\cdot,y)\|_\infty E_{m-1}(f^*)\right. \\ & \left.+ \| h(\cdot,f^*)\|_\infty \, \sup_{|y| \leq 1} E_{m-1}(k_2(\cdot,y))  + 2 \sup_{|y| \leq 1} \|k_2(\cdot,y)\|_\infty E_{m-1}(h(\cdot,f^*))\right].
\end{align*}
Now, by the assumptions  \eqref{assK} on the kernels $k_i$ and taking \eqref{erroreSobolev} into account, we get
\begin{align*}
\|f^* -f_m\|_\infty &\leq  \mathcal{C} \sup_{|y| \leq 1} \|k_1(\cdot,y)\|_{W^r} \left[ \frac{\|f^*\|_{\infty}}{m^r}+2  E_{m-1}(f^*)   \right] \\  &+ \sup_{|y| \leq 1} \|k_2(\cdot,y)\|_{W^r} \left[ \frac{1}{m^r} \| h(\cdot,f^*)\|_\infty    + 2  E_{m-1}(h(\cdot,f^*))\right] .
\end{align*}
Moreover, note that by the hypothesis \eqref{assgh} on $g$ and by virtue of Theorem \ref{teo:K}, we can deduce that $f^* \in W^r$. Therefore, we can apply Theorem \ref{teo:BestApprx}. Then, we obtain
\begin{align*}
\|f^* -f_m\|_\infty  
& \leq  \frac{\mathcal{C}}{m^r} \left( \sup_{|y| \leq 1} \|k_1(\cdot,y)\|_{W^r} \|f^*\|_{W^r} \right. \\ & \left.+ \sup_{|y| \leq 1} \|k_2(\cdot,y)\|_{W^r} \left[ \| h(\cdot,f^*)\|_\infty    + 2^r B_r  \|h\|_{\mathbf{W}^r(\Omega)} \|f^*\|^s_{W^r}\right] \right), 
\end{align*}
from which the assertion.
\end{proof}
\begin{proof}\emph{of Theorem {\ref{teo:conv*}}.}  
	The proof can be leaded following step by step the proof of Theorem \ref{teo:conv}, just by substituting operator $K_2$ and $K_m^2$ with $K^\mu$ defined in (\ref{Kmu}) and $K^*_m$ defined in (\ref{defK*m}), respectively, and arriving to the following inequality, also using (\ref{prodotto}),
\begin{eqnarray}	
\|f^* -f_m\|_\infty &\leq  &\mathcal{C}\left[ \|({K}^1- {K}^1_m)f^*\|_\infty+  \|(\widetilde{K}^\mu-{K}_m^*)f^*\|_\infty \right] \nonumber
	\\ & =& \mathcal{C} \left[ \sup_{|y| \leq 1} |e_m(k_1(\cdot,y) f^*)|+\sup_{|y| \leq 1} |r_m (h(\cdot,f^*),y)| \right] \nonumber \\ & \leq &\C \left[\sup_{|y| \leq 1} E_{2m-1}(k_1(\cdot,y) f^*)+
	\sup_{|y| \leq 1} E_{m-1}( h(\cdot,f^*)) \right] \label{f-fm1}
	\end{eqnarray}
	where $\mathcal{C}\neq{C}(m,f)$.
	
	For the first term in the brackets, under the assumptions \eqref{assK1}, using (\ref{jackson}), and the same arguments in the proof of Theorem \ref{teo:conv}, we get 
\begin{equation}\label{Ek1}
\sup_{|y| \leq 1} E_{2m-1}(k_1(\cdot,y) f^*)\leq \mathcal{C}\frac{\|f^*\|_{Z^s}}{m^s}.
\end{equation}
	Concerning the second term, by \eqref{jackson} we have
	\begin{equation}\label{Eh}
	E_{m-1}( h(\cdot,f^*)) \leq \mathcal{C} \int_0^{\frac{1}{m}} \frac{\Omega^k_\varphi(h(\cdot,f^*),t)}{t} dt,
	\end{equation}
	and then it is crucial to estimate the modulus of smoothness $\Omega_\varphi^k(h(\cdot,f^*),t)$. Under the assumptions we made, we can consider $k=1$. 
	By the definition \eqref{Omega}, we have
	\begin{align}\label{T1+T2}
&	|\Delta_{\tau \varphi} h(x,f^*(x))| \nonumber \\&=\left|h\left(x+\tau \frac{\varphi(x)}{2}, f^*\left(x+\tau \frac{\varphi(x)}{2} \right)\right)-h\left(x-\tau \frac{\varphi(x)}{2}, f^*\left(x-\tau \frac{\varphi(x)}{2} 
	 \right)\right)\right| \nonumber \\ &\leq \left|h\left(x+\tau \frac{\varphi(x)}{2}, f^*\left(x+\tau \frac{\varphi(x)}{2} \right)\right)-h\left(x-\tau \frac{\varphi(x)}{2}, f^*\left(x+\tau \frac{\varphi(x)}{2} 
	 \right)\right)\right| \nonumber \\ &+ \left|h\left(x-\tau \frac{\varphi(x)}{2}, f^*\left(x+\tau \frac{\varphi(x)}{2} \right)\right)-h\left(x-\tau \frac{\varphi(x)}{2}, f^*\left(x-\tau \frac{\varphi(x)}{2} 
	 \right)\right) \right| \nonumber \\ &=T_1(x)+T_2(x).
	\end{align}
By using the derivability of $h$ with respect to the first variable, we get
	\begin{equation}\label{T1}
	T_1(x) \leq \tau \left|\frac{\partial h}{\partial x}\left(\xi_1, f^*\left(x+\frac{\tau}{2} \varphi(x)\right)\right) \right|, \qquad \xi_1 \in \left[x-\frac{\tau}{2}\varphi(x), x+\frac{\tau}{2}\varphi(x) \right] 
	\end{equation}
and similarly, by the derivability of $h$ with respect to the 	second variable we write
	\begin{equation}\label{T2}
	T_2(x) \leq \left|\frac{\partial h}{\partial y}\left(x-\frac{\tau}{2}\varphi(x), f^*(\xi_2)\right) \right| \left|\Delta_{\tau \varphi} f^*(x) \right|, \,  \xi_2 \in \left[x-\frac{\tau}{2}\varphi(x), x+\frac{\tau}{2}\varphi(x) \right]. 
	\end{equation}
	Hence, by combining \eqref{T1} and \eqref{T2} in \eqref{T1+T2} and considering that we are assuming  $h \in W^1$, we have
	\begin{equation*} 
	|\Delta_{\tau \varphi} h(x,f^*(x))| \leq \|h\|_{W^1} [\tau+ \left|\Delta_{\tau \varphi} f^*(x) \right|],
	\end{equation*}
	from which we can conclude that 
$$\Omega_\varphi(h(\cdot,f^*),t) \leq \|h\|_{W^1} \left[ t+ \Omega_\varphi(f,t) \right].$$
Therefore, since by Theorem \ref{classeDebSing}, it is $f^* \in Z^s$, by \eqref{Eh} we get
$$E_{m-1}(h(\cdot,f^*)) \leq \mathcal{C} \|h\|_{W^1} \left[\frac{1}{m}+ \frac{1}{m^s} \right] \leq  \frac{\C}{m^s} \|h\|_{W^1}.$$
The assertion follows by combining the above result and \eqref{Ek1} into \eqref{f-fm1}.
\end{proof}


\section{Declarations}
The authors are member of the INdAM Research group GNCS and the TAA-UMI Research Group. This research has been accomplished within ``Research ITalian network on Approximation'' (RITA). 

L. Fermo, C. Laurita, and M.G. Russo are partially supported by the PRIN 2022 PNRR project no. P20229RMLB
financed by the European Union - NextGeneration EU and by the Italian Ministry of University and Research (MUR). 

L. Fermo is also partially supported by the INdAM-GNCS project 2024 ``Algebra lineare numerica per problemi di grandi dimensioni: aspetti teorici e applicazioni''.

The authors have no other competing interests to declare that are relevant to the content of this article.

%
%


\begin{thebibliography}{10}
\providecommand{\url}[1]{{#1}}
\providecommand{\urlprefix}{URL }
\expandafter\ifx\csname urlstyle\endcsname\relax
  \providecommand{\doi}[1]{DOI~\discretionary{}{}{}#1}\else
  \providecommand{\doi}{DOI~\discretionary{}{}{}\begingroup
  \urlstyle{rm}\Url}\fi

\bibitem{KEAtkinson}
Atkinson, K.E.: The numerical evaluation of fixed points for completely
  continuous operators.
\newblock SIAM J. Numer. Anal. \textbf{10}, 799--807 (1973)

\bibitem{Atkinsonsurvey}
Atkinson, K.E.: A survey of numerical methods for solving nonlinear integral
  equations.
\newblock J. Integral Equations Appl. \textbf{4}, 1--32 (1992)

\bibitem{A}
Atkinson, K.E.: The Numerical Solution of Integral Equations of the Second
  Kind.
\newblock Cambridge University Press, Cambridge (1997)

\bibitem{AtkChan}
Atkinson, K.E., Chandler, G.: Boundary integral equation methods for solving
  {L}aplce's equation with nonlinear boundary conditions: the smooth boundary
  case.
\newblock Math. Comp. \textbf{55}(192), 451--472 (1990)

\bibitem{barrera22}
Barrera, D., Bartoň, M., Chiarella, I., Remogna, S.: On numerical solution of
  {F}redholm and {H}ammerstein integral equations via {N}ystr\"om method and
  {G}aussian quadrature rules for splines.
\newblock Appl. Numer. Math. \textbf{174}, 71--88 (2022)

\bibitem{barrera20}
Barrera, D., El~Mokhtari, F., Ibáñez, M.J., Sbibih, D.: Non-uniform
  quasi-interpolation for solving {H}ammerstein integral equations.
\newblock Int. J. Comput. Math. \textbf{97}(1--2), 72--84 (2020)

\bibitem{bialecki1981}
Bialecki, R., Nowak, A.J.: Boundary value problems in heat conduction with
  nonlinear material and nonlinear boundary conditions.
\newblock Appl. Math. Model. \textbf{5}(6), 417--421 (1981)

\bibitem{Chidume}
Chidume, C., Adamu, A., Okereke, L.: Iterative algorithms for solutions of
  {H}ammerstein equations in real {B}anach spaces.
\newblock Fixed Point Theory Appl. \textbf{12} (2020)

\bibitem{dagnino19}
Dagnino, C., Dallefrate, A., Remogna, S.: Spline quasi-interpolating projectors
  for the solution of nonlinear integral equations.
\newblock J. Comput. Appl. Math. \textbf{354}, 360--372 (2019)

\bibitem{das2016}
Das, P., Sahani, M.M., Nelakanti, G., Long, G.: Legendre spectral projection
  methods for {F}redholm--{H}ammerstein integral equations.
\newblock J. Sci. Comput. \textbf{68}, 213--230 (2016)

\bibitem{debonis06}
De~Bonis, M., Mastroianni, G.: Projection methods and condition numbers in
  uniform norm for {F}redholm and {C}auchy singular integral equations.
\newblock SIAM J. Numer. Anal. \textbf{44}(4), 1351--1374 (2006)

\bibitem{DT}
Ditzian, Z., Totik, V.: Moduli of smoothness.
\newblock Springer-Verlag, New York (1987)

\bibitem{Dolezale}
Dolezale, V.: Monotone Operators and Its Applications in Automation and Network
  Theory.
\newblock Studies in Automation and Control. Elsevier, New York (1979)

\bibitem{ElJeSlSt}
Elschner, J., Jeon, Y., Sloan, I.H., Stephan, E.P.: The collocation method for
  mixed boundary value problem on domains with curved polygonal boundaries.
\newblock Numer. Math. \textbf{76}, 355--381 (1997)

\bibitem{FermoLaurita2015}
Fermo, L., Laurita, C.: A {N}ystr\"om method for a boundary integral equation
  related to the {D}irichlet problem on domains with corners.
\newblock Numer. Math. \textbf{130}(1), 35--71 (2015)

\bibitem{FermoLaurita}
Fermo, L., Laurita, C.: On the numerical solution of a boundary integral
  equation for the exterior {N}eumann problem on domains with corners.
\newblock Appl. Numer. Math. \textbf{94}, 179--200 (2015)

\bibitem{FermoLaurita2}
Fermo, L., Laurita, C.: A {N}ystr\"om method for mixed boundary value problems
  in domains with corners.
\newblock Appl. Numer. Math. \textbf{149}, 65--82 (2020)

\bibitem{BestApproxErr}
Fermo, L., Laurita, C., Russo, M.: On the error of best polynomial
  approximation of composite functions.
\newblock submitted  (2023)

\bibitem{FermoRusso2010}
Fermo, L., Russo, M.: Numerical methods for {F}redholm integral equations with
  singular right-hand sides.
\newblock Adv. Comput. Math. \textbf{33}, 305--330 (2010)

\bibitem{Kaneko90}
Kaneko, H., Noren, R., Xu, Y.: Regularity of the solution of {H}ammerstein
  equations with weakly singular kernel.
\newblock Integral Equations Operator Theory \textbf{13}(5), 660 – 670 (1990)

\bibitem{Kaneko91}
Kaneko, H., Xu, Y.: Degenerate kernel method for {H}ammerstein equations.
\newblock Math. Comp. \textbf{56}(193), 141--148 (1991)

\bibitem{KantarovicAkilov}
Kantorovich, L., Akilov, G.: Functional Analysis.
\newblock Pergamon Press, Oxford (1982)

\bibitem{kelmanson1984}
Kelmanson, M.A.: Solution of nonlinear elliptic equations with boundary
  singularities by an integral equation method.
\newblock J. Comput. Phys. \textbf{56}(2), 244--258 (1984)

\bibitem{Krasnoseski64}
Krasnosel'skii, M.A.: Topological Methods in the Theory of Nonlinear Integral
  Equations.
\newblock Macmillan, New York (1964)

\bibitem{Krasnoseski}
Krasnosel'skii, M.A., Zabreiko, P.: Geometrical Methods of Nonlinear Analysis.
\newblock Springer-Verlag, Berlin (1984)

\bibitem{Kress}
Kress, R.: Linear Integral Equations, \emph{Applied Mathematical Sciences},
  vol.~82.
\newblock Springer-Verlag, Berlin (1989)

\bibitem{Kumar90}
Kumar, S.: The numerical solution of {H}ammerstein equations by a method based
  on polynomial collocation.
\newblock J. of Austral. Math. Soc. Ser. B. \textbf{31}(3), 319–329 (1990)

\bibitem{KumarSloan}
Kumar, S., Sloan, I.H.: A new collocation-type method for {H}ammerstein
  integral equations.
\newblock Math. Comp. \textbf{48}, 585--593 (1987)

\bibitem{Laurita}
Laurita, C.: A numerical method for the solution of exterior {N}eumann problems
  for the {L}aplace equation in domains with corners.
\newblock Appl. Numer. Math. \textbf{119}, 248--270 (2017)

\bibitem{Laurita2020}
Laurita, C.: A new stable numerical method for {M}ellin integral equations in
  weighted spaces with uniform norm.
\newblock Calcolo \textbf{57}(3) (2020)

\bibitem{LR}
Luther, U., Russo, M.: Boundedness of the {H}ilbert transformation in some
  weighted {B}esov type spaces.
\newblock Integral Equations Operator Theory \textbf{36}(2), 220--240 (2000)

\bibitem{Mandal21}
Mandal, M., Kant, K., Nelakanti, G.: Discrete {L}egendre spectral methods for
  {H}ammerstein type weakly singular nonlinear {F}redholm integral equations.
\newblock Int. J. Comput. Math. \textbf{98}(11), 2251--2267 (2021)

\bibitem{Mandal1}
Mandal, M., Nelakanti, G.: Superconvergence results for weakly singular
  {F}redholm-{H}ammerstein integral equations.
\newblock Numer. Funct. Anal. Optim. \textbf{40}(5), 548--570 (2019)

\bibitem{MastroianniMilovanovic}
Mastroianni, G., Milovanovic, G.: Interpolation Processes: Basic Theory and
  Applications.
\newblock Springer Monographs in Mathematics (2008)

\bibitem{MRT}
Mastroianni, G., Russo, M., Themistoclakis, W.: Numerical {M}ethods for
  {C}auchy {S}ingular {I}ntegral {E}quations in {S}paces of {W}eighted
  {C}ontinuous {F}unctions.
\newblock In: I.~Gohberg, al. (eds.) Recent Advances in Operator Theory and its
  Applications, pp. 311--336. Birkh{\"a}user Basel, Basel (2005)

\bibitem{Mezzanotte2021}
Mezzanotte, D., Occorsio, D., Russo, M.G.: Combining {N}ystr\"om methods for a
  fast solution of {F}redholm integral equations of the second kind.
\newblock Mathematics \textbf{9}(21), 2652 (2021)

\bibitem{MS}
Monegato, G., Scuderi, L.: Global polynomial approximation for {S}ymm's
  equations on polygons.
\newblock Numer. Math. \textbf{86}, 655--683 (2000)

\bibitem{MonegatoScuderi}
Monegato, G., Scuderi, L.: A polynomial collocation method for the numerical
  solution of weakly singular and singular integral equations on non-smooth
  boundaries.
\newblock Int. J. Numer. Meth. Engng \textbf{58}, 1985--2011 (2003)

\bibitem{Narendra}
Narendra, K., Gallman, P.: An iterative method for the identification of
  nonlinear systems using a {H}ammerstein model.
\newblock IEEE Trans. Autom. Control \textbf{11}(3), 546--550 (1966)

\bibitem{Pascali}
Pascali, D., Sburlan, S.: Nonlinear mappings of monotone type.
\newblock Springer (1978)

\bibitem{PV}
Pedas, A., Vainikko, G.: Superconvergence of piecewise polynomial collocations
  for nonlinear weakly singular integral equations.
\newblock Journal of Integral Equations and Applications \textbf{9}(4),
  379--406 (1997)

\bibitem{Pot2022}
P\"otzsche, C.: Urysohn and {H}ammerstein operators on {H}\"older spaces.
\newblock Analysis \textbf{42}(4), 205--240 (2022)

\bibitem{RKP}
Rashidinia, J., Khosravian~Arab, H., Parsa, A.: Gauss-{L}egendre quadrature for
  solution of {H}ammerstein integral equations.
\newblock Math. Sci. \textbf{5}(4), 345--354 (2011)

\bibitem{ruotsalainen1988}
Ruotsalainen, K., Wendland, W.: On the boundary element method for some
  nonlinear boundary value problems.
\newblock Numer. Math. \textbf{53}, 299--314 (1988)

\bibitem{Scuderi}
Scuderi, L.: A {C}hebyshev polynomial collocation {BIEM} for mixed boundary
  value problems on nonsmooth boundaries.
\newblock J. Integral Equations Appl. \textbf{14}, 179--221 (2002)

\bibitem{Shav}
Shahsavaran, A., Fotros, F.: An effective and simple scheme for solving
  nonlinear {F}redholm integral equations.
\newblock Math. Model. Anal. \textbf{27}(2), 215--231 (2022)

\bibitem{szego}
Szeg\"o, G.: Orthogonal Polynomials, \emph{Amer. Math. Soc.}, vol.~23.
\newblock Amer. Math. Soc. Colloq. Publ., Providence (1975)

\bibitem{Wazwaz}
Wazwaz, A.: Linear and Nonlinear Integral Equations: Methods and Applications,
  1st edn.
\newblock Springer Publishing Company, Incorporated (2011)

\end{thebibliography}

%
%

\end{document}